\documentclass[final,a4paper,leqno,11pt]{article}
\usepackage{graphicx}
\usepackage{mathptmx}
\usepackage{amssymb,amsmath, amsfonts,amsthm}
\usepackage{a4wide}
\usepackage{color}

\newcommand{\email}[1]{{\tt #1}}
\newcommand{\R}{\mathbb{R}}
\newcommand{\N}{\mathbb{N}}
\newcommand{\norm}[1]{\|#1\|}

\newcommand{\dist}[1]{{\rm dist}(#1)}
\newcommand{\B}{{\cal B}}

\newcommand{\K}{{\cal K}}

\newcommand{\Q}{{\cal Q}}
\newcommand{\HH}{{\cal H}}

\newcommand{\setto}[1]{\mathop{\rightarrow}\limits^#1}

\newcommand{\skalp}[1]{\langle #1\rangle}
\newcommand{\bskalp}[1]{\big\langle #1\big\rangle}

\newcommand{\xb}{\bar x}
\newcommand{\yb}{\bar y}
\newcommand{\zb}{\bar z}
\newcommand{\ub}{\bar u}

\newcommand{\lb}{\bar\lambda}

\newcommand{\xba}{{\bar x^\ast}}

\newcommand{\oo}{o}
\newcommand{\OO}{{O}}
\newcommand{\argmin}{{\rm arg\,min\,}}
\newcommand{\argmax}{\mathop{\rm arg\,max}\limits}

\newcommand{\inn}{{\rm int\,}}
\newcommand{\bd}{{\rm bd\,}}
\newcommand{\gph}{\mathrm{gph}\,}
\newcommand{\tto}{\rightrightarrows}
\def\h{\hfill\Box}
\def\disp{\displaystyle}
\def\dn{\downarrow}
\def\la{\langle}
\def\ra{\rangle}
\def\st{\stackrel}
\def\O{\Omega}

\def\Hat{\widehat}
\def\ox{\bar{x}}

\def\oy{\bar{y}}

\def\ov{\bar{v}}

\def\gph{\mbox{\rm gph}\,}
\def\dist{\mbox{\rm dist}\,}
\def\epi{\mbox{\rm epi}\,}

\def\dom{\mbox{\rm dom}\,}
\def\ker{\mbox{\rm ker}\,}

\def\lip{\mbox{\rm lip}\,}

\def\bd{\mbox{\rm bd}\,}

\def\oR{\Bar{\R}}
\def\lm{\lambda}
\def\gg{\gamma}
\def\dd{\delta}
\def\al{\alpha}

\def\kk{\kappa}

\def\tilde{\widetilde}

\def\ph{\varphi}
\def\hat{\Hat}

\def\tto{\;{\lower 1pt \hbox{$\rightarrow$}}\kern -10pt
\hbox{\raise 2pt \hbox{$\rightarrow$}}\;}
\def\Hat{\widehat}

\def\Bar{\overline}
\def\ra{\rangle}
\def\la{\langle}

\def\epsilon{\varepsilon}
\def\B{\Bbb B}
\def\h{\hfill\Box}
\def\R{\Bbb R}

\def\N{\Bbb N}
\def\ox{\bar{x}}

\def\oy{\bar{y}}

\def\ov{\bar{v}}

\def\gph{\mbox{\rm gph}\,}

\def\epi{\mbox{\rm epi}\,}

\def\dom{\mbox{\rm dom}\,}

\def\ker{\mbox{\rm ker}\,}

\def\lip{\mbox{\rm lip}\,}

\def\h{\hfill\square}
\def\dn{\downarrow}
\def\O{\Omega}
\def\ph{\varphi}
\def\emp{\emptyset}
\def\st{\stackrel}
\def\oR{\Bar{\R}}
\def\lm{\lambda}
\def\gg{\gamma}

\def\dd{\delta}
\def\al{\alpha}

\newlength{\myparboxwidth}
\newtheorem{theorem}{Theorem}[section]

\newtheorem{lemma}[theorem]{Lemma}

\newtheorem{definition}[theorem]{Definition}
\newtheorem{example}[theorem]{Example}

\title{Characterizations of Tilt-Stable Minimizers\\
in Second-Order Cone Programming}
\author{MAT\'U\v{S} BENKO\footnote{Institute of Computational Mathematics, Johannes Kepler University Linz, A-4040 Linz, Austria; \email{matusbenko@hotmail.com}}\and HELMUT GFRERER\footnote{Institute of Computational Mathematics, Johannes Kepler University Linz, A-4040 Linz, Austria; \email{helmut.gfrerer@jku.at}}
$\;$ \and $\;$ BORIS S. MORDUKHOVICH\footnote{Department of Mathematics, Wayne State University, Detroit, MI 48202, USA;
\email{boris@math.wayne.edu}}}
\date{}
\begin{document}
\maketitle
\small
{\bf Abstract.} This paper is devoted to the study of tilt stability of local minimizers, which plays an important role in both theoretical and numerical aspects of optimization. This notion has been comprehensively investigated in the unconstrained framework as well as for problems of nonlinear programming with $C^2$-smooth data. Available results for nonpolyhedral conic programs were obtained only under strong constraint nondegeneracy assumptions. Here we develop an approach of second-order variational analysis, which allows us to establish complete neighborhood and pointbased characterizations of tilt stability for problems of second-order cone programming generated by the nonpolyhedral second-order/Lorentz/ice-cream cone. These characterizations are established under the weakest metric subregularity constraint qualification condition.
\vspace*{0.05in}

\vspace*{0.05in}

\vspace*{0.05in}

{\bf Abbreviated title.} Tilt stability in conic programming\vspace*{-0.15in}

\normalsize
\section{Introduction and Brief Overview}\vspace*{-0.05in}

This paper mainly aims at deriving efficient characterizations of tilt-stable local minimizers in problems of {\em second-order cone programming} (SOCP) given as follows:
\begin{equation}\label{EqProbl}
\mbox{minimize }\;f(x)\;\mbox{ subject to }\;x\in\Gamma:=\big\{x\in\R^n\;\big|\;g(x)\in\Q\big\},
\end{equation}
where $f\colon\R^n\to\R$ and $g\colon\R^n\to\R^{1+m}$ are $C^2$-smooth around the reference point $\xb$, and where
\begin{equation}\label{2nd-cone}
\Q:=\big\{(q_0,q_r)\in\R\times\R^{m}\;\big|\;\norm{q_r}-q_0\le 0\big\}
\end{equation}
is known as the {\em second-order/Lorentz/ice-cream cone}. Problems of this type, often abbreviated as SOCPs, have been well recognized in conic programming and numerous applications; see, e.g., \cite{ag03,BonSh00} for more discussions and references. Observe that, despite the imposed smoothness of $f$ and $g$, the defined SOCP model is a problem of nonsmooth optimization due to the nondifferentiability at the origin of the norm function presented in (\ref{2nd-cone}).

The notion of tilt stability of local minimizers was introduced by Poliquin and Rockafellar \cite{PolRo98} in the extended-real-valued framework of unconstrained optimization as follows. Given a function $\ph\colon\R^n\to\oR:=(-\infty,\infty]$ with the domain $\dom\ph:=\{x\in\R^n\;|\;\ph(x)<\infty\}$, we say that $\ox\in\dom\ph$ is a {\em tilt-stable minimizer} of $\ph$ if there is a number $\gg>0$ such that the mapping
\begin{equation}\label{EqM_gamma}
M_\gamma(v^\ast):=\argmin\big\{\ph(x)-\skalp{v^\ast,x}\;\big|\;\|x-\ox\|\le\gg\big\},\quad v^*\in\R^n,
\end{equation}
is single-valued and Lipschitz continuous in some neighborhood of $\ov^*:=0\in\R^n$ with $M_\gamma(0)=\{\xb\}$. The main result of \cite{PolRo98} provides a characterization of tilt-stable local minimizers for $\ph$ via the positive-definiteness of the second-order subdifferential/generalized Hessian of $\ph$ at the reference point as introduced by the third author in \cite{Mo92}. The results of \cite{PolRo98} were extended by Mordukhovich and Nghia \cite{mor-nghia13,MoNg14} who developed a new approach to tilt stability and derived {\em quantitative} characterization of tilt-stable local minimizers with a prescribed Lipschitz modulus of the mapping $M_\gg$ in \eqref{EqM_gamma} and then computed the exact lower bound of these moduli in terms of the second-order subdifferential and its ``combined" modification in finite and infinite dimensions. The characterizations of \cite{mor-nghia13,MoNg14} were established in both {\em neighborhood} form involving points nearby the reference minimizer as well as via {\em pointbased} criteria expressed entirely at the point in question.

Quite recently, Chieu et al. \cite{ChHiNg17} derived quantitative neighborhood characterizations of tilt stability in the abstract finite-dimensional framework of \eqref{EqM_gamma} that are expressed in the form of \cite{MoNg14} but with replacing the dual-space combined second-order subdifferential therein by the primal-dual ``subgradient graphical derivative" discussed in Section~2. We refer the reader to \cite{BonSh00,DrLe13,DMN14,m18,mor-nghia13,MoNg14,ZN14} and the bibliographies therein for other results on tilt stability as well as closely related notions in the abstract framework of \eqref{EqM_gamma}. Applications of tilt stability to multiplier criticality and convergence rates of some primal-dual algorithms in numerical optimization have been recently provided in \cite{ms16}.

Efficient implementations and applications of the second-order characterizations of tilt stability obtained in the abstract scheme \eqref{EqM_gamma} to structural classes of optimization problems with explicit constraints require developing adequate machinery of second-order variational analysis and subdifferential calculus under appropriate qualification conditions. First results in this direction were obtained by Mordukhovich and Rockafellar \cite{MoRo12} for problems of nonlinear programming (NLPs) with $C^2$-smooth data under the classical {\em linear independence constraint qualification} (LICQ) fulfilled at the given local minimizer. In this setting, as shown in \cite{MoRo12}, tilt stability is characterized by Robinson's {\em strong second-order sufficient condition} (SSOSC) introduced in \cite{Rob80}.

Since LICQ is not necessary for tilt stability in nonlinear programming (in contrast to the case of strong regularity of the associated KKT systems), several attempts have been made to relax LICQ and study tilt-stable local minimizers for NLPs with nonunique Lagrange multipliers; see \cite{ChHiNg17,gm15,MoNg14,MoOut13}. The most advanced results in this direction were obtained under the weakest {\em metric subregularity constraint qualification} (MSCQ). It was first used by Gfrerer and Mordukhovich \cite{gm15}, being combined with some other weak constraint qualifications, to derive pointbased sufficient conditions, necessary conditions as well as complete characterizations of tilt-stable minimizers for NLPs. More recently, Chieu at al. \cite{ChHiNg17} employed MSCQ to establish a new neighborhood characterization and pointbased sufficient conditions for tilt stability in NLPs with $C^2$-smooth inequality constraints.\vspace*{0.03in}

All the methods and results for tilt stability in NLPs discussed above are strongly based on the polyhedral structure of such problems. Not much is known for tilt stability of local minimizers in {\em nonpolyhedral} problems of conic programming. The main result of \cite{mos} yields a characterization of tilt-stable minimizers for SOCPs in the extended form of SSOSC under a {\em nondegeneracy} condition corresponding to LICQ in the SOCP setting. It is based on the second-order subdifferential calculations taken from Outrata and Sun \cite{os08}. Similar characterizations of tilt stability involving nondegeneracy and uniqueness of Lagrange multipliers are established in \cite{mnr15,mor15} for $C^2$-reducible problems of conic programming, where the second-order subdifferential term is calculated in \cite{mnr15}, based on the results of \cite{dsy14}, for semidefinite programs entirely in terms of the program data.\vspace*{0.03in}

The main goal of this paper is to derive verifiable sufficient conditions, necessary conditions, and complete characterizations of tilt-stable minimizers for SOCPs under an appropriate version of {\em MSCQ} in second-order cone programming that is introduced here. Our results include neighborhood characterizations and much more difficult {\em pointbased} conditions for tilt stability. Furthermore, we establish pointbased {\em quantitative} evaluations and precise formulas for computing the {\em exact bound} of tilt stability expressed entirely in terms of the given SOCP data. All the results obtained are the first in the literature for tilt stability in second-order cone programming without imposing nondegeneracy. In particular, they are new under the Robinson constraint qualification, which corresponds to the replacement of metric subregularity by the stronger metric regularity assumption. Our results provide new information (e.g., the calculation of the exact bound of tilt stability) even under nondegeneracy.\vspace*{0.03in}

The rest of the paper is organizes as follows. Section~2 presents some preliminaries from variational analysis and generalized differentiation systematically used in the subsequent text. In Section~3 we obtain neighborhood characterizations of tilt-stable minimizers for SOCPs with neighborhood calculating the exact bound of tilt stability. In Section~4 we begin developing significantly more challenging and more important pointbased results on tilt stability of local minimizers in SOCPs starting with second-order sufficient conditions for such minimizers under MSCQ. The major example presented here illustrates essential features of the obtained conditions in comparison with those known under nondegeneracy and polyhedrality. Section~5 is devoted to deriving ``no-gap" pointbased necessary conditions for tilt-stable minimizers in SOCPs and their complete characterizations for which an additional ``2-regularity" assumption is required in some cases. In the concluding Section~6 we discuss the main thrust of the paper and some open problems in this area.\vspace*{0.03in}

Our notation is conventional in variational analysis, conic programming, and generalized differentiation; see, e.g., \cite{BonSh00,m18,RoWe98}. Recall that $\B$ and ${\cal S}$ stand respectively for the closed unit ball and sphere of the space in question, that $\B_r(x):=x+r\B$, and that $\N:=\{1,2,\ldots\}$. Taking into account the structure of $\Q$ in \eqref{2nd-cone}, we represent an element $q\in\R^{1+m}$ as $q=(q_0,q_r)$ with
$q_0\in\R$ and $q_r\in\R^m$, and also denote $\hat q:=(-q_0,q_r)$. This notation is used, in particular, for the mapping $g(x)=\big(g_0(x),g_r(x)\big)$ in the constraint system $\Gamma$ from \eqref{EqProbl}. Furthermore, we have the following representations for the dual cone, the tangent cone, and the normal cone to $\Q$, respectively: $\Q^*=\{\hat q\mid q\in\Q\}$,
\begin{subequations}\label{EqTangNormalConeQ}
\begin{gather}
T_\Q(q)=\begin{cases}\Q&\text{if $q=0$},\\
\R^{1+m}&\text{if $q\in\inn\Q$},\\
\Big\{u\in\R^{1+m}\;\Big|\;\disp\frac{q_r}{\norm{q_r}}u_r-u_0\le 0\Big\}&\text{if $q\in\bd\Q\setminus\{0\}$},
\end{cases}\\
N_\Q(q)=\begin{cases}\Q^*&\text{if $q=0$},\\
\{0\}&\text{if $q\in\inn\Q$},\\
\Big\{\alpha\Big(-1,\disp\frac{q_r}{\norm{q_r}}\Big)\;\Big|\;\alpha\ge 0\Big\}=\big\{\alpha\hat q\;\big|\;\alpha\ge 0\big\}&\text{if $q\in\bd\Q\setminus\{0\}$}.
\end{cases}\end{gather}
\end{subequations}
Given a (sufficiently) smooth real-valued/scalar function $\ph:\R^n\to\R$, denote its gradient and Hessian at $x$ by $\nabla\ph(x)$ and $\nabla^2\ph(x)$, respectively. Considering further a vector function $h\colon\R^n\to\R^s$ with $s>1$, denote by $\nabla h(x)$ the Jacobian of $h$ at $x$, while by $\nabla^2 h(x)$ we mean a linear mapping from $\R^n$ into the space of $s\times n$ matrices defined by
$$
\nabla^2 h(x)u:=\disp\lim_{t\to 0}\frac{\nabla h(x+tu)-\nabla h(x)}t.
$$
We also employ the second-order construction $\nabla^2h(x)(u,v):=(\nabla^2 h(x)u)v=\lim_{t\to 0}\frac{\nabla h(x+tu)-\nabla h(x)}tv$.\vspace*{-0.15in}

\section{Preliminaries from Variational Analysis}\vspace*{-0.05in}

We first recall generalized differential constructions of variational analysis utilized in what follows. The reader can find more details in the books \cite{m18,RoWe98} and the references therein. Starting with sets, consider a nonempty one $\O\subset\R^n$ and define the (Bouligand--Severi) {\em tangent cone} to $\O$ at $\ox\in\O$ by
\begin{equation}\label{tan}
T_\O(\ox):=\{v\in\R^n\;\big|\;\exists\,t_k\dn 0,\;\exists\,v_k\to v\;\mbox{ such that }\;\ox+t_k v_k\in\O,\;k\in\N\big\}.
\end{equation}
The (limiting, Mordukhovich) {\em normal cone} to $\O$ at $\ox\in\O$ is given by
\begin{equation}\label{nor}
N_\O(\ox):=\Big\{x^*\in\R^n\;\Big|\;\exists\,x_k\st{\O}\to\ox,\;\exists\,x^*_k\to x^*\;\mbox{ such that }\;\disp\limsup_{x\to\ox}\frac{\la x^*,x-x_k\ra}{\|x-x_k\|}\le 0\Big\},
\end{equation}
where $x\st{\O}\to\ox$ indicates that $x\to\ox$ with $x\in\O$. Note that both cones \eqref{tan} and \eqref{nor} may be nonconvex.

Considering next an extended-real-valued function $\ph\colon\R^n\to\oR$ finite at $\ox$, the only {\em subdifferential} of $\ph$ at $\ox$ used in this paper is defined by
\begin{equation}\label{sub}
\partial\ph(\ox):=\Big\{x^*\in\R^n\;\Big|\;\exists\,x_k\st{\ph}\to\ox,\;\exists\,x^*_k\to x^*\;\mbox{ such that }\;\disp\liminf_{x\to x_k}\frac{\ph(x)-\ph(x_k)-\la x^*_k,x-x_k\ra}{\|x-x_k\|}\ge 0\Big\},
\end{equation}
where $x\st{\ph}\to\ox$ indicates that $x\to\ox$ with $\ph(x)\to\ph(\ox)$. Note that despite the nonconvexity of $\partial\ph(\ox)$, the subdifferential \eqref{sub} enjoys full calculus based on variational/extremal principles of variational analysis \cite{m18,RoWe98}. There are close relationships between \eqref{sub} and \eqref{nor}, namely:
\begin{equation}\label{rel}
\partial\ph(\ox)=\big\{x^*\in\R^n\;\big|\;(x^*,-1)\in N_{{\rm\small epi}\,\ph}\big(\ox,\ph(\ox)\big)\big\}\;\mbox{ and }\;N_\O(\ox)=\partial\dd_\O(\ox),
\end{equation}
where $\epi\ph:=\{(x,\al)\in\R^{n+1}\;|\;\al\ge\ph(x)\}$ is the epigraph of $\ph$, and where $\dd_\O(x)$ is the indicator function of $\ph$ equal to 0 if $x\in\O$ and to $\infty$ otherwise.

Recall that a lower semicontinuous (l.s.c.) function $\ph\colon\R^n\to\oR$ is {\em continuously prox-regular} at $\ox$ for $\ox^*\in\partial\ph(\ox)$ if there exists $r>0$ such that
\begin{eqnarray*}
\ph(x)\ge\ph(u)+\la x^*,x-u\ra-(r/2)\|x-u\|^2\;\mbox{ whenever }\;x,u\;\mbox{ near }\;\ox\;\mbox{ and }\;x^*\in\partial\ph(u)\;\mbox{ near }\;\ox^*
\end{eqnarray*}
and if the mapping $(x,x^*)\mapsto\ph(x)$ is continuous relative to the subdifferential graph $\gph\partial\ph$ at $(\ox,\ox^*)$. This class is rather large including \cite{RoWe98} all l.s.c.\ convex functions, strongly amenable functions, etc.

Considering further a set-valued mapping $F\colon\R^n\tto\R^m$, recall that the {\em graphical derivative} of $F$ at $(\ox,\ov)\in\gph F$ is defined via the tangent cone \eqref{tan} by
\begin{equation}\label{der}
DF(\ox,\ov)(u):=\big\{w\in\R^m\;\big|\;(u,w)\in T_{{\rm\small gph}\,F}(\ox,\ov)\big\},\quad u\in\R^n.
\end{equation}
The main {\em second-order} construction for extended-real-valued functions used in this paper is the following {\em primal-dual} one \cite{mor15} known as the {\em subgradient graphical derivative} of $\ph\colon\R^n\to\oR$ at $\ox\in\dom\ph$ for $\ox^*\in\partial\ph(\ox)$, which is defined via \eqref{sub} and \eqref{der} by
\begin{equation}\label{2der}
D\partial\ph(\ox,\ox^*)(u):=\big\{w\in\R^n\;\big|\;(u,w)\in T_{{\rm\small gph}\,\partial\ph}(\ox,\ox^*)\big\},\quad u\in\R^n.
\end{equation}

The generalized second-order derivative \eqref{2der} is exploited in \cite[Theorem 3.3]{ChHiNg17} to derive the following neighborhood characterization of tilt-stable minimizers in the abstract framework of extended-real-valued functions. As in \cite{mor-nghia13}, the {\em exact bound of tilt stability} for $\ph$ at $\ox$ is defined by
\begin{equation*}
{\rm tilt}(\ph,\ox)=\inf_{\gg>0}\lip M_\gg(0),
\end{equation*}
where $\lip M_\gg(0)$ stands for the infimum of Lipschitz moduli of $M_\gg$ from \eqref{EqM_gamma} around the origin.\vspace*{-0.05in}

\begin{theorem}[\bf abstract neighborhood characterization of tilt stability]\label{ThNeighbCharTiltStab}
Let $\ph\colon\R^n\to\oR$ be a l.s.c.\ function that is continuously prox-regular at $\xb\in\dom\ph$ for $\xba=0$. Assume that $0\in\partial\ph(\ox)$. Then $\xb$ is a tilt-stable local minimizer of $\ph$ with modulus $\kappa>0$ if and only if there is a constant $\eta>0$ such that
\begin{equation}\label{abs}
\skalp{u^*,u}\ge\frac 1\kappa\text{ whenever }\;u\in{\cal S},\;u^*\in D\partial\ph(x,x^*)(u),\;(x,x^*)\in\gph\partial\ph\cap\B_\eta(\xb,0).
\end{equation}
Furthermore, the exact bound for tilt stability of $\ph$ at $\ox$ is calculated by the formula
\begin{equation}\label{exact}
{\rm tilt}(\ph,\xb)=\inf_{\eta>0}\left\{\frac{1}{\skalp{u^*,u}}\;\Big|\;u\in{\cal S},\;u^*\in D\partial\ph(x,x^*)(u),\;(x,x^*)\in\gph\partial\ph\cap\B_\eta(\xb,0)\right\}.
\end{equation}
\end{theorem}\vspace*{-0.05in}

At what follows we aim at employing this result to establish verifiable conditions for and characterizations of tilt stability in SOCPs with expressing them entirely via the {\em program data}. The main attention is paid to deriving {\em pointbased conditions} formulated exactly at the reference point. Our major assumptions is the novel {\em metric subregularity constraint qualification} that is far removed from constraint nondegeneracy being essentially weaker than the conventional Robinson constraint qualification under which the obtained results are also new. The developed techniques are rather involved and are based on second-order calculus and calculations in nonpolyhedral settings.

Note that to study tilt stability in the SOCP framework, it is sufficient considering only the case where $g(\ox)=0$, which we {\em always assume in what follows} without further mentioning. Indeed, otherwise it holds that either $g(\xb)\in\inn\Q$, or $g(\xb)\in\bd\Q\setminus\{0\}$. In the first case the constraint in \eqref{EqProbl} can be ignored in our local analysis. In the second case it can be equivalently written as a $C^2$-inequality $\norm{g_{r}(x)}-g_{0}(x)\le 0$ for which there exists a comprehensive tilt stability theory.\vspace*{-0.15in}

\section{Neighborhood Characterizations of Tilt-Stable Minimizers for SOCPs}\vspace*{-0.05in}

We start with neighborhood characterizations of tilt stability for SOCPs and first introduce the underlying metric subregularity constraint qualification for the constraint system in \eqref{EqProbl}. Recall that a mapping $Q\colon\R^n\tto\R^m$ is {\em metrically regular} around $(\ox,\oy)\in\gph Q$ with modulus $\sigma>0$ if there exist neighborhoods $U$ of $\ox$ and $V$ of $\oy$ such that we have the estimate
\begin{equation}\label{mr}
{\rm dist}\big(x;Q^{-1}(y)\big)\le\sigma\,{\rm dist}\big(y;Q(x)\big)\;\mbox{ for all }\;x\in U\;\mbox{ and }\;y\in V,
\end{equation}
where dist$(w;\O)$ signifies the distance between a point and a set. The mapping $F$ is {\em metrically subregular} at $(\ox,\oy)$ with modulus $\sigma$ if the distance estimate \eqref{mr} holds with $y=\oy$ therein.\vspace*{-0.05in}

\begin{definition}[\bf metric subregularity constraint qualification for SOCPs]\label{mscq} We say that the {\sc metric subregularity constraint qualification} $(MSCQ)$ is fulfilled for the SOCP constraint system $\Gamma$ from \eqref{EqProbl} at a point $\ox\in\Gamma$ with modulus $\sigma>0$ if the mapping $Q(x):=g(x)-\Q$ is metrically subregular at $(\ox,0)$ with this modulus, i.e.,
\begin{equation}\label{subreg}
{\rm dist}\big(x;g^{-1}(\Q)\big)\le\sigma\,{\rm dist}\big(g(x);\Q\big)\;\mbox{ for all }\;x\in U.
\end{equation}
\end{definition}

To this end, it is worth emphasizing that the more restrictive metric {\em regularity} of the mapping $Q(\cdot)$ from Definition~\ref{mscq} around $(\ox,0)$ is equivalent to the {\em Robinson constraint qualification}
\begin{equation*}
0\in{\rm int}\big\{g(\ox)+\nabla g(\ox)\R^n-\Q\big\}.
\end{equation*}

Taking into account that the validity of MSCQ at $\ox\in\Gamma$ yields its fulfillment at any $x\in\Gamma$ sufficiently close to $\ox$ (with the same modulus $\sigma$ as can be supposed without loss of generality), we deduce from
\cite[Theorem~4.1]{HenJouOut02}, \cite[Proposition~1]{HenOut05}, and the formulas in \eqref{EqTangNormalConeQ} that the tangent cone and the normal cone to $\Gamma$ at $x$ are computed by
\begin{align}
\nonumber T_\Gamma(x)&=\big\{u\in\R^n\;\big|\;\nabla g(x)u\in T_\Q\big(g(x)\big)\big\}\\
\label{EqTanGamma}&=\begin{cases}\big\{u\in\R^n\;\big|\;\nabla g(x)u\in\Q\}&\text{if $g(x)=0$},\\
\R^n&\text{if $g(x)\in\inn\Q$},\\
\Big\{u\in\R^n\;\Big|\;\disp\Big\la\frac{g_r(x)}{\norm{g_r(x)}},\nabla g_r(x)u\Big\ra-\nabla g_0(x)u\le 0\Big\}&\text{if $g(x)\in\bd\Q\setminus\{0\}$};
\end{cases}
\end{align}
\[
N_\Gamma(x)=\begin{cases}\big\{\nabla g(x)^*\lambda\;\big|\;\lambda\in\Q^*\big\}&\text{if $g(x)=0$},\\
\{0\}&\text{if $g(x)\in\inn\Q$},\\
\left\{\alpha\Big(\nabla g_r(x)^*\disp\frac{g_r(x)}{\norm{g_r(x)}}-\nabla g_0(x)^*\Big)\;\Big|\;\alpha\ge 0\right\}&\text{if $g(x)\in\bd\Q\setminus\{0\}$}.
\end{cases}
\]

Given a normal vector $x^*\in N_\Gamma(x)$, consider the corresponding {\em multiplier set}
$$
\Lambda(x,x^\ast):=\big\{\lambda\in N_\Q(g(x))\;\big|\;\nabla g(x)^*\lambda=x^*\big\}
$$
and define the {\em critical cone} to $\Gamma$ at $(x,x^*)$ by
$$
\K_\Gamma(x,x^*)=T_\Gamma(x)\cap[x^*]^\perp.
$$
Note that for every $\lambda\in\Lambda(x,x^*)$ we have $u\in\K_\Gamma(x,x^*)$ if and only if $u\in T_\Gamma(x)$ and $\skalp{\lambda,\nabla g(x)u}=0$.
If $x^\ast=0$, the only multiplier important for our analysis is $\lambda=0$. Thus we introduce the set
$$
\Lambda^0(x,x^*):=\begin{cases}\{0\}&\text{when $x^*=0$},\\
\Lambda(x,x^*)&\text{otherwise.}\end{cases}
$$
Further, for every direction $u\in\R^n$ define the {\em directional multiplier set} by
$$
\Lambda(x,x^*;u):=\argmax_{\lambda\in\lambda(x,x^*)}\bskalp{\big(\nabla^2\skalp{\lambda, g}(x)+\HH(x,\lambda)\big)u,u},
$$
where the {\em curvature mapping} $\HH\colon\R^n\times\R^{1+m}\to\R^{n\times n}$ is given by
\begin{equation}\label{H}
\HH(x,\lambda):=\begin{cases}\disp\frac{-\lambda_0}{g_0(x)}\Big(\nabla g_r(x)^*\nabla g_r(x)-\nabla g_0(x)^*\nabla g_0(x)\Big)&\text{if $g(x)\in\bd \Q\setminus\{0\}$},\\
0&\text{otherwise.}\end{cases}
\end{equation}

Now we are ready to establish neighborhood characterizations of tilt stability for local minimizers of SOCPs with computing the exact bound of tilt stability via the given data of \eqref{EqProbl} and \eqref{2nd-cone}. Tilt stability of $\ox$ in \eqref{EqProbl} is naturally understood as the one for $\ph(x):=f(x)+\delta_\Gamma(x)$ in the sense defined above. Then $\ox\in\Gamma$ is a {\em stationary point} of \eqref{EqProbl} if $0\in\partial\ph(\ox)$, i.e.,
\begin{equation*}
0\in\partial f(\ox)+N_{\Gamma}(\ox)
\end{equation*}
due to the standard first-order subdifferential sum rule.\vspace*{-0.05in}

\begin{theorem}[\bf neighborhood characterizations of tilt stability for SOCPs]\label{ThNeighCharRedSyst} Let $\ox\in\Gamma$ be a stationary point of \eqref{EqProbl}, and let MSCQ hold at $\ox$ with modulus $\sigma$ from \eqref{subreg}. Then the following three assertions are equivalent:

{\bf(i)} The point $\xb$ is tilt-stable minimizer for \eqref{EqProbl} with modulus $\kappa>0$.

{\bf(ii)} There exists a constant $\eta>0$ such that
\begin{equation}\label{EqNeighborTilt}
\left\{\parbox{13cm}{$\bskalp{\big(\nabla^2f(x)+\nabla^2\skalp{\lambda,g}(x)+\HH(x,\lambda)\big)u,u}\ge\disp\frac 1\kappa$
\\whenever $x\in\Gamma\cap\B_\eta(\xb)$, $\norm{x^\ast}\le\eta$, $u\in\K_\Gamma\big(x,x^\ast-\nabla f(x)\big)\cap{\cal S}$,\\
$\lambda\in\Lambda\big(x,x^\ast-\nabla f(x);u\big)\cap\sigma\norm{x^\ast-\nabla f(x)}\B$.}\right.
\end{equation}

{\bf(iii)} There exists a constant $\eta>0$ such that
\begin{equation}\label{EqNeighborTilt1}
\left\{\parbox{13cm}{$\bskalp{\big(\nabla^2f(x)+\nabla^2\skalp{\lambda,g}(x)+\HH(x,\lambda)\big)u,u}\ge\disp\frac 1\kappa$
\\whenever $x\in\Gamma\cap\B_\eta(\xb)$, $\norm{x^\ast}\le\eta$, $u\in\K_\Gamma\big(x,x^\ast-\nabla f(x)\big)\cap{\cal S}$,\\
\mbox{\qquad}$\lambda\in\Lambda\big(x,x^\ast-\nabla f(x);u\big)\cap\Lambda^0\Big(x,x^\ast-\nabla f(x)\big)$.}\right.
\end{equation}
Furthermore, the exact bound for tilt stability of $\ox$ in \eqref{EqProbl} is calculated by
\begin{equation}\label{exact0}
\begin{array}{ll}
{\rm tilt}(f+\delta_\Gamma,\xb)=\disp\inf_{\eta>0}\sup\Bigl\{1\big/\bskalp{\big(\nabla^2f(x)+\nabla^2\skalp{\lambda, g}(x)+\HH(x,\lambda)\big)u,u}\;\Big|\;x\in\Gamma\cap\B_\eta(\xb),\\
\norm{x^\ast}\le\eta,\;u\in\K_\Gamma\big(x,x^\ast-\nabla f(x)\big)\cap{\cal S},\;\lambda\in\Lambda\big(x,x^\ast-\nabla f(x);u\big)\cap\sigma\norm{x^\ast-\nabla f(x)}\B\Bigr\}.
\end{array}
\end{equation}
\end{theorem}\vspace*{-0.05in}
{\bf Proof.} We derive the claimed results from those in Theorem~\ref{ThNeighbCharTiltStab} with $\ph(x)=f(x)+\delta_\Gamma(x)$ by using the appropriate rules of first-order and second-order generalized differential calculus. Since $f\in C^2$, it follows that  $\partial(f+\delta_\Gamma)(x)=\nabla f(x)+N_\Gamma(x)$ and
$$
D\partial(f+\delta_\Gamma)(x,x^*)=\nabla^2f(x)+DN_\Gamma\big(x,x^*-\nabla f(x)\big)
$$
for all $x$ sufficiently close to $\xb$. Due to the validity of MSCQ at such $x$ with modulus $\sigma$, we get from the second-order calculation in \cite[Corollary~5.2]{HaMoSa17} that
\begin{align*}
DN_\Gamma\big(x,x^*-\nabla f(x)\big)(u)=&\Bigl\{\nabla^2\skalp{\lambda,g}(x)u +\HH(x,\lambda)u\;\big|\;
\lambda\in\Lambda\big(x,x^\ast-\nabla f(x);u\big)\cap\sigma\norm{x^*-\nabla f(x)}\B\Bigr\}\\
&\qquad+N_{\K_\Gamma(x,x^*-\nabla f(x))}(u)\;\mbox{ for any }\;u\in\R^n.
\end{align*}
Noting that the validity of MSCQ at $x\in\Gamma$ with modulus $\sigma$ ensures that this property holds at $x$ with any modulus $\sigma'\ge\sigma$ gives us the formula
\begin{align*}
DN_\Gamma\big(x,x^*-\nabla f(x)\big)(u)=&\Bigl\{\nabla^2\skalp{\lambda,g}(x)u +\HH(x,\lambda)u\;\Big|\;\lambda\in\Lambda\big(x,x^\ast-\nabla f(x);u\big)\cap\Lambda^0\big(x,x^\ast-\nabla f(x)\big)\B\Bigr\}\\
&\qquad+N_{\K_\Gamma(x,x^*-\nabla f(x))}(u)\;\mbox{ for any }\;u\in\R^n.
\end{align*}
Substituting the above calculations into \eqref{abs} and \eqref{exact} with taking into account that for every vector $\zeta\in N_{\K_\Gamma(x,x^*-\nabla f(x))}(u)$ we have $\skalp{\zeta,u}=0$ allows us to complete the proof of the theorem. $\h$\vspace*{-0.15in}

\section{Pointbased Sufficient Conditions for Tilt Stability in SOCP}\vspace*{-0.05in}

The results of Theorem~\ref{ThNeighCharRedSyst} involve {\em all} the points from a neighborhood of the reference stationary point $\ox$, which makes the verification of the obtained formula \eqref{EqNeighborTilt} difficult. The goal of this section is to derive verifiable sufficient conditions for tilt stability in SOCPs expressed only at the point $\ox$ in question. We first present the the following two technical lemmas. \vspace*{-0.05in}

\begin{lemma}[\bf curvature mapping on critical directions]\label{LemHH} For every pair $(x,x^*)\in\gph N_\Gamma$, every multiplier $\lambda\in\Lambda(x,x^*)$, and every critical direction $u\in \K_\Gamma(x,x^*)$ we have
\[
\bskalp{\HH(x,\lambda)u,u}\ge 0.
\]
\end{lemma}\vspace*{-0.05in}
{\bf Proof.} It can be assumed without loss of generality that $\lambda_0\not=0$ and that $g(x)\in\bd\Q\setminus\{0\}$, i.e., $g_0(x)=\norm{g_r(x)}>0$. Indeed, otherwise  the lemma conclusion holds trivially. Since $\lambda\in N_\Q(g(x))$, it follows that $\lambda_0=-\norm{\lambda_r}<0$. Further, we deduce from $u\in\K_\Gamma(x,x^*)$ and $\lambda\in\Lambda(x,x^*)$ that
\begin{align*}
0&=\skalp{\lambda,\nabla g(x)u}=\lambda_0\nabla g_0(x)u+\skalp{\lambda_r,\nabla g_r(x)u}=-\norm{\lambda_r}\nabla g_0(x)u+\skalp{\lambda_r,\nabla g_r(x)u}\\
&\begin{cases}\le\norm{\lambda_r}\big(-\nabla g_0(x)u+\norm{\nabla g_r(x)u}\big)&\\
\ge\norm{\lambda_r}\big(-\nabla g_0(x)u-\norm{\nabla g_r(x)u}\big)&
\end{cases}
\end{align*}
implying therefore that $(\nabla g_0(x)u)^2\le\norm{\nabla g_r(x)u}^2$. Using now the identity
$$
\bskalp{\HH(x,\lambda)u,u}=\disp\frac{-\lambda_0}{g_0(x)}\Big(\norm{\nabla g_r(x)u}^2-(\nabla g_0(x)u)^2\Big)
$$
leads us to the claimed inequality and thus completes the proof. $\h$\vspace*{-0.05in}

\begin{lemma}[\bf limiting procedure]\label{LemSequ} Consider convergent sequences $x_k\setto{\Gamma}\xb$, $x_k^*\to 0$, $\lambda^k\to\tilde\lambda$, and $u_k\to\tilde u$ as $k\to\infty$ satisfying the inclusions
\[
\lambda^k\in\Lambda\big(x_k,x_k^*-\nabla f(x_k);u_k\big)\cap\sigma\norm{x_k^*-\nabla f(x_k)}\B\;\mbox{ and }\;u_k\in\K_\Gamma\big(x_k,x_k^*-\nabla f(x_k)\big)\cap{\cal S}
\]
for all $k\in\N$ together with the condition
\begin{equation}\label{EqBoundedH}
\limsup_{k\to\infty}\skalp{\HH(x_k,\lambda^k)u_k,u_k}<\infty.
\end{equation}
Then we have the limiting relationships
\begin{subequations}\label{EqPropTilde_u_lambda}
\begin{align}&\tilde\lambda\in\Lambda\big(\xb,-\nabla f(\xb)\big)\cap\sigma\norm{\nabla f(\xb)}\B\subset\Lambda^0\big(\xb,-\nabla f(\xb)\big),\\
&\tilde u\in{\cal S},\;\skalp{\tilde\lambda,\nabla g(\xb)\tilde u}=0,\ \tilde\lambda_0\big(\norm{\nabla g_r(\xb)\tilde u}^2-(\nabla g_0(\xb)\tilde u)^2\big)=0.
\end{align}
\end{subequations}
\end{lemma}\vspace*{-0.05in}
{\bf Proof.} Passing to the limit as $k\to\infty$ in the conditions
$$
\nabla g(x_k)^*\lambda^k=x_k^\ast-\nabla f(x_k),\;\lambda^k\in N_\Q\big(g(x_k)\big)\subset\Q^*,\;\mbox{ and }\;\norm{\lambda^k}\le\sigma\norm{x_k^*-\nabla f(x_k)}
$$
verifies that $\nabla g(\xb)^*\tilde\lambda^k=-\nabla f(\xb)$, $\tilde\lambda\in\Q^*$, and $\norm{\tilde\lambda}\le\sigma\norm{\nabla f(\xb)}$. This yields
$$
\tilde\lambda\in\Lambda\big(\xb,-\nabla f(\xb)\big)\cap\sigma\norm{\nabla f(\xb)}\B\subset\Lambda^0\big(\xb,-\nabla f(\xb)\big).
$$
It follows from $\skalp{\lambda^k,\nabla g(x_k)u_k}=0$ as $k\in\N$  that $\skalp{\lb,\nabla g(\xb)\tilde u}=0$. Thus it remains to
show that
\begin{equation}\label{lim0}
\tilde\lambda_0\big(\norm{\nabla g_r(\xb)\tilde u}^2-(\nabla g_0(\xb)\tilde u)^2\big)=0.
\end{equation}
Since \eqref{lim0} is certainly true if $\tilde\lambda_0=0$, we proceed with the case where $\tilde\lambda_0\not=0$ and hence $\lambda^k_0\not=0$ for all large $k\in\N$. It follows from $\tilde\lambda\in\Q^*$ that $-\tilde\lambda_0\ge\norm{\tilde\lambda_r}$. If $g(x_k)=0$ for infinitely many $k$, then $\nabla g(x_k)u_k\in\Q$ and therefore $\nabla g(\xb)\tilde u\in\Q$, which in turn is equivalent to $\nabla g_0(\xb)\tilde u\ge\norm{\nabla g_r(\xb)\tilde u}$. From $\skalp{\tilde\lambda,\nabla g(\xb)\tilde u}=0$ we deduce the conditions
\[
-\tilde\lambda_0\nabla g_0(\xb)\tilde u=\skalp{\tilde\lambda_r,\nabla g_r(\xb)\tilde u}\le\norm{\tilde\lambda_r}\cdot\norm{\nabla g_r(\xb)\tilde u}\le-\tilde\lambda_0\nabla g_0(\xb)\tilde u
\]
implying that $\norm{\tilde\lambda_r}\cdot\norm{\nabla g_r(\xb)\tilde u}=-\tilde\lambda_0\nabla g_0(\xb)\tilde u$. Thus we get that either $\norm{\tilde\lambda_r}=-\tilde\lambda_0$ and $\norm{\nabla g_r(\xb)\tilde u}=\nabla g_0(\xb)\tilde u$, or $\norm{\nabla g_r(\xb)\tilde u}=\nabla g_0(\xb)\tilde u=0$, which verifies equality \eqref{lim0} in this case. In the remaining case where $g(x_k)\not=0$ for all but finitely many $k\in\N$ it follows that $g(x_k)\in\bd\Q\setminus\{0\}$ for such $k$. This is due to the fact that the condition $0\not=\lambda^k\in N_\Q(g(x_k))$ yields
\[
\skalp{\HH(x_k,\lambda^k)u_k,u_k}=\frac{-\lambda_0^k}{g_0(x_k)}\Big(\norm{\nabla g_r(x_k)u_k}^2-(\nabla g_0(x_k)u_k)^2\Big).
\]
Using \eqref{EqBoundedH} together with $g_0(x_k)\to 0$ as $k\to\infty$ and Lemma~\ref{LemHH} tells us that
\[
\lim_{k\to\infty}\lambda_0^k\big(\norm{\nabla g_r(x_k)u_k}^2-(\nabla g_0(x_k)u_k)^2\big)=\tilde\lambda_0\big(\norm{\nabla g_r(\xb)\tilde u}^2-(\nabla g_0(\xb)\tilde u)^2\big)=0,
\]
which verifies \eqref{lim0} in the latter case and thus completes the proof. $\h$\vspace*{0.05in}

Now we are ready to proceed with deriving pointbased sufficient conditions for tilt stability in SOCP. It is instructive to split our consideration into the two cases: {\bf (a)} {\em in-kernel} $\K_\Gamma(\xb,-\nabla f(\xb))\subset\ker\nabla g(\xb)$ and {\bf(b)} {\em out-of-kernel} $\K_\Gamma(\xb,-\nabla f(\xb))\not\subset\ker\nabla g(\xb)$. The second case can be treated directly by passing to the limit in the conditions of Theorem~\ref{ThNeighCharRedSyst} with the usage of the lemmas above.\vspace*{-0.05in}

\begin{theorem}[\bf sufficient condition for tilt stability in the out-of-kernel case]\label{ThSuffCondNonZeroCrit} In addition to the assumptions of Theorem~{\rm\ref{ThNeighCharRedSyst}} suppose that there is $\ub\in\K_\Gamma(\xb,-\nabla f(\xb))$ satisfying $\nabla g(\xb)\ub\not=0$. Denote
\begin{equation}\label{barl}
\lb:=\frac{\norm{\nabla f(\xb)}}{\norm{\nabla g(\xb)^*\nabla\hat g(\xb)\ub}}\nabla\hat g(\xb)\ub\quad\mbox{ with }\quad\hat g(x):=\big(-g_0(x),g_r(x)\big)
\end{equation}
and assume also that for some $\kappa>0$ we have
\begin{equation}\label{EqPointTilt1}
\begin{array}{ll}
\bskalp{\big(\nabla^2 f(\xb)+\nabla^2\skalp{\lb,g}(\xb)\big)u,u}>\disp\frac{1}{\kappa}\quad\mbox{whenever}\\
u\in{\cal S}\;\mbox{ with }\;\skalp{\lb,\nabla g(\xb)u}=0,\;\lb_0\big(\norm{\nabla g_r(\xb)u}^2-(\nabla g_0(\xb)u)^2\big)=0.
\end{array}
\end{equation}
Then $\xb$ is a tilt-stable local minimizer for \eqref{EqProbl} with modulus $\kappa$. Moreover, we have the upper estimate
\begin{align*}
{\rm tilt}(f+\delta_\Gamma,\xb)\le\sup\Bigl\{1\big/\bskalp{\big(\nabla^2f(\xb)+\nabla^2\skalp{\lb,g}(\xb)\big)u,u}\;\Big|\;&u\in{\cal S},\;\skalp{\lb,\nabla g(\xb)u}=0,\\
&\lb_0\big(\norm{\nabla g_r(\xb)u}^2-(\nabla g_0(\xb)u)^2\big)=0\Bigr\}.
\end{align*}
\end{theorem}\vspace*{-0.05in}
{\bf Proof.} Employing Theorem~\ref{ThNeighCharRedSyst}, it suffices to show that the second-order condition \eqref{EqPointTilt1} yields \eqref{EqNeighborTilt}. Then the claimed upper estimate of the exact bound of tilt stability follows easily. Assuming on the contrary that \eqref{EqNeighborTilt} fails while \eqref{EqPointTilt1} holds, we find sequences $x_k\setto{\Gamma}\xb$, $x_k^*\to 0$, $u_k\in \K_\Gamma(x_k,x_k^*-\nabla f(x_k))\cap{\cal S}$, and $\lambda^k\in\Lambda(x_k,x_k^*-\nabla f(x_k);u_k)\cap\sigma\norm{x_k^*-\nabla f(x_k)}\B$ satisfying
\begin{equation}\label{EqAux1}
\bskalp{\big(\nabla^2f(x_k)+\nabla^2\skalp{\lambda^k,g}(x_k)+\HH(x_k,\lambda^k)\big)u_k,u_k}<\frac 1\kappa\;\mbox{ for all}\;k\in\N,
\end{equation}
which ensures, in particular, that the sequence of $\skalp{\HH(x_k,\lambda^k)u_k,u_k}$ is bounded. Thus passing to a subsequence if necessary tells us that $u_k$ converges to some $\tilde u\in{\cal S}$ and that $\lambda_k$ converges to some $\tilde\lambda$ satisfying \eqref{EqPropTilde_u_lambda} by Lemma~\ref{LemSequ}. Employing further
$$
\Lambda\big(\xb,-\nabla f(\xb)\big)\subset\Q^*\cap\big(\nabla g(\xb)\ub\big)^\perp=N_\Q\big(\nabla g(\xb)\ub\big)\subset\big\{\alpha\nabla\hat g(\xb)\ub\;\big|\;\alpha\ge 0\big\}
$$
together with $\nabla g(\xb)^*\lambda=-\nabla f(\xb)$ for all $\lambda\in\Lambda(\xb,-\nabla f(\xb))$ ensures that $\Lambda^0(\xb,-\nabla f(\xb))=\{\lb\}$ and hence $\tilde\lambda=\lb$. It follows from \eqref{EqAux1} and Lemma~\ref{LemHH} that
\[
\bskalp{\big(\nabla^2 f(x)+\nabla^2\skalp{\lb,g}(x)\big)\tilde u,\tilde u}\le \limsup_{k\to\infty}\bskalp{\big(\nabla^2f(x_k)+\nabla^2\skalp{\lambda^k, g}(x_k)+\HH(x_k,\lambda^k)\big)u_k,u_k}\le\frac 1\kappa,
\]
which contradicts \eqref{EqPointTilt1} and thus completes the proof of the theorem. $\h$\vspace*{0.05in}

For our subsequent analysis in the in-kernel case we introduce the set
\begin{equation}\label{z}
\begin{array}{ll}
{\cal Z}:=\Big\{(u,\lambda,v,w)\in{\cal S}\times\Lambda^0\big(\xb,-\nabla f(\xb)\big)\times{\cal S}\times\R^n\;\Big|\;\skalp{\lambda,\nabla g(\xb)u}=0,\\
\lambda_0\big(\norm{\nabla g_r(\xb)u}^2-(\nabla g_0(\xb)u)^2\big)=0,\;\nabla g(\xb)v=0,\;\lambda\in N_{\Q}\big(\nabla g(\xb)w+\frac 12\nabla^2 g(\xb)(v,v)\big)\Big\}
\end{array}
\end{equation}
Further, for every triple $(u,\lambda,v)\in\R^n\times\Q^*\times\R^n$ define the {\em infimum function}
\begin{align}\label{EqRho}
\rho(u,\lambda,v):=\inf_{z\in\R^n}\Big\{&-\lambda_0\Big(\norm{\nabla g_r(\xb)z+\nabla^2g_r(\xb)(v,u)}^2-\big(\nabla g_0(\xb)z+\nabla^2g_0(\xb)(v,u)\big)^2\Big)\;\Big|\\
\nonumber&\skalp{\lambda,\nabla g(\xb)z+\nabla^2g(\xb)(v,u)}=0\Big\}.
\end{align}
with the convention that $\rho(u,\lambda,v):=\infty$ whenever $\{z\in\R^n\;|\;\skalp{\lambda,\nabla g(\xb)z+\nabla^2g(\xb)(v,u)}=0\}=\emp$.\vspace*{0.02in}

Let us present another lemma before deriving the main result of this section in the in-kernel case.\vspace*{-0.05in}

\begin{lemma}[\bf properties of the infimum function]\label{LemRho} For any triple $(u,\lambda,v)\in\R^n\times\Q^*\times\R^n$ we get that $\rho(u,\lambda,v)\ge0$. Furthermore, the infimum in \eqref{EqRho} is attained whenever $\rho(u,\lambda,v)$ is finite.
\end{lemma}\vspace*{-0.05in}
{\bf Proof.} If $\lambda=0$, we clearly have that $\rho(u,\lambda,v)=0$ and that the infimum in \eqref{EqRho} is attained at any $z\in\R^n$. Considering now the case where $\lambda\not=0$ and $\rho(u,\lambda,v)<\infty$, pick any $z$ satisfying $\skalp{\lambda,\nabla g(\xb)z+\nabla^2g(\xb)(v,u)}=0$ and denote $\eta:=\nabla g(\xb)z+\nabla^2g(\xb)(v,u)$. It follows from $\lambda\in\Q^*\setminus\{0\}$ and the structure of $\Q$ in \eqref{2nd-cone} that $-\lambda_0\ge\norm{\lambda_r}$ and $-\lambda_0>0$. Hence
\[
-\lambda_0\vert\eta_0\vert=\vert\skalp{\lambda_r,\eta_r}\vert\le\norm{\lambda_r}\cdot\norm{\eta_r}\le-\lambda_0\norm{\eta_r},
\]
which implies in turn the relationships
$$
0\le-\lambda_0(\norm{\eta_r}^2-\eta_0^2)=-\lambda_0\Big(\norm{\nabla g_r(\xb)z+\nabla^2g_r(\xb)(v,u)}^2-\big(\nabla g_0(\xb)z+\nabla^2g_0(\xb)(v,u)\big)^2\Big),
$$
and therefore $\rho(u,\lambda,v)\ge 0$. To show finally that the infimum is attained in \eqref{EqRho}, observe that the minimization problem therein is with a quadratic cost and one linear equality constraint. It can be equivalently reduced by some linear transformation $z=b+Ay$ to an unconstrained quadratic program of the form $\min\frac 12\la y,By\ra+\skalp{c,y}$; see, e.g., \cite[Chapter~10.1]{Fle81}. The boundedness of the cost function ensures that the matrix $B$ is positive semidefinite and that there is some vector $\yb$ satisfying the first-order optimality condition $B\yb=-c$. Indeed, assuming the opposite yields by the fundamental theorem of linear algebra the existence of some direction $d$ with $Bd=0$ and $\skalp{-c,d}\ne 0$, which contradicts the boundedness of the cost function. Thus the unconstrained optimization problem admits a solution $\yb$, and so $\zb:=b+A\yb$ is an optimal solution to the constrained problem under consideration. $\h$\vspace*{0.05in}

The next theorem is a major result of the paper that establishes pointbased sufficient conditions for tilt-stable minimizers for SOCPs in the most involved in-kernel case.\vspace*{-0.05in}

\begin{theorem}[\bf sufficient conditions for tilt stability in the in-kernel case]\label{ThSuffCondZeroCrit} In addition to the standing assumptions of Theorem~{\rm\ref{ThNeighCharRedSyst}} suppose that
\begin{equation}\label{EqCritDir0}
\nabla g(\xb)u=0\;\mbox{ for all }\;u\in\K_\Gamma\big(\xb,-\nabla f(\xb)\big).
\end{equation}
Assume also that there is a number $\kappa>0$ such that the following two conditions are satisfied:

{\bf (a)} For every $u\in K_\Gamma(\xb,-\nabla f(\xb))\cap{\cal S}$ and every $\lambda\in\Lambda(\xb,-\nabla f(\xb);u)$ we have
\[
\bskalp{\big(\nabla^2f(\xb)+\nabla^2\skalp{\lambda,g}(\xb)\big)u,u}>\frac 1\kappa.
\]

{\bf(b)} For every quadruple $(u,\lambda,v,w)\in{\cal Z}$ we have
\begin{equation}\label{EqSecOrder1}
\begin{array}{ll}
\disp\bskalp{\big(\nabla^2f(\xb)+\nabla^2\skalp{\lambda,g}(\xb)\big)u,u}+\frac{\rho(u,\lambda,v)}{\nabla g_0(\xb)w+\frac 12\nabla^2 g_0(\xb)(v,v)}>\frac 1\kappa\\
\disp\mbox{ whenever }\;\nabla g(\xb)w+\frac 12\nabla^2 g(\xb)(v,v)\not=0;\quad\mbox{ and}
\end{array}
\end{equation}\vspace*{-0.15in}
\begin{equation}\label{EqSecOrder2}
\disp\bskalp{\big(\nabla^2f(\xb)+\nabla^2\skalp{\lambda,g}(\xb)\big)u,u}>\frac 1\kappa\;\mbox{ whenever }\;\rho(u,\lambda,v)=0.
\end{equation}
Then $\xb$ is a tilt-stable local minimizer for \eqref{EqProbl} with modulus $\kappa$. Moreover, we have the upper estimate
\begin{equation}\label{tilt upper}
{\rm tilt}(f+\delta_\Gamma,\xb)\le\frac 1{\min\{\chi_1,\chi_2,\chi_3\}},
\end{equation}
where the numbers $\chi_i,\;i=1,2,3$, are calculated by
\begin{subequations}
\begin{align*}
\chi_1&:=\inf\big\{\bskalp{\big(\nabla^2f(\xb)+\nabla^2\skalp{\lambda,g}(\xb)\big)u,u}\;\big|\;u\in K_\Gamma\big(\xb,-\nabla f(\xb)\big)\cap{\cal S},\; \lambda\in\Lambda\big(\xb,-\nabla f(\xb);u\big)\big\},\\
\chi_2&:=\inf\big\{\bskalp{\big(\nabla^2f(\xb)+\nabla^2\skalp{\lambda,g}(\xb)\big)u,u}\;\big|\;\exists v,w\;\mbox{ with }\;(u,\lambda,v,w)\in{\cal Z},\;\rho(u,\lambda,v)=0\big\},\\
\chi_3&:=\inf\Big\{\bskalp{\big(\nabla^2f(\xb)+\nabla^2\skalp{\lambda,g}(\xb)\big)u,u}+\frac{\rho(u,\lambda,v)}{\nabla g_0(\xb)w+\frac 12\nabla^2 g_0(\xb)(v,v)}\;\Big|\\
\nonumber&\qquad\quad(u,\lambda,v,w)\in{\cal Z},\;\nabla g(\xb)w+\frac 12\nabla^2 g(\xb)(v,v)\not=0\Big\}
\end{align*}
with the set ${\cal Z}$ defined in \eqref{z}.
\end{subequations}
\end{theorem}\vspace*{-0.05in}
{\bf Proof.} Let us show that the fulfillment of both conditions (a) and (b) implies that the tilt stability characterization \eqref{EqNeighborTilt} of Theorem~\ref{ThNeighCharRedSyst} holds. This verifies all the conclusions of the theorem including the upper bound estimate \eqref{tilt upper}, which follows from \eqref{exact0} and the proof below.

Arguing by contradiction, suppose that \eqref{EqNeighborTilt} fails, i.e., there are sequences $x_k\setto{\Gamma}\xb$, $x_k^*\to 0$, $u_k\in\K_\Gamma(x_k,x_k^*-\nabla f(x_k))\cap{\cal S}$, and $\lambda^k\in\Lambda(x_k,x_k^*-\nabla f(x_k);u_k)\cap\sigma\norm{x_k^*-\nabla f(x_k)}\B$ satisfying
\begin{equation}\label{EqAux2}
\bskalp{\big(\nabla^2f(x_k)+\nabla^2\skalp{\lambda^k, g}(x_k)+\HH(x_k,\lambda^k)\big)u_k,u_k}<\frac 1\kappa\;\mbox{ for all }\;k\in\N,
\end{equation}
and then show that either condition (a) or condition (b) is violated. Passing to a subsequence if necessary, we get that the sequences $u_k$ and $\lambda^k$ converge to some $\tilde u$ and $\tilde\lambda$, respectively. Since \eqref{EqAux2} yields the boundedness of $\skalp{\HH(x_k,\lambda^k)\big)u_k,u_k}$, it follows from Lemma~\ref{LemSequ} that
\[
\tilde u\in{\cal S},\;\tilde\lambda\in\Lambda^0\big(\xb,-\nabla f(\xb)\big),\;\skalp{\tilde\lambda,\nabla g(\xb)\tilde u}=0,\;\tilde\lambda_0\big(\norm{\nabla g_r(\xb)\tilde u}^2-(\nabla g_0(\xb)\tilde u)^2\big)=0.
\]
Let us split the subsequent proof into the following two cases:\\[1ex]
{\bf Case~I:} {\em We have  that $x^k=\xb$ for infinitely many $k$}. Suppose without loss of generality that it holds for all $k\in\N$. Then $\nabla g(\xb)u^k\in\Q$ whenever $k\in\N$, and thus $\nabla g(\xb)\tilde u\in\Q$ yielding $\tilde u\in T_\Gamma(\xb)$ by \eqref{EqTanGamma} due to the imposed MSCQ. Together with $\skalp{\tilde\lambda,\nabla g(\xb)\tilde u}=0$ it tells us that $\tilde u\in\K_\Gamma(\xb,-\nabla f(\xb))$. Next we show that $\tilde\lambda\in\Lambda(\xb,-\nabla f(\xb);\tilde u)$. Assuming on the contrary that $\tilde\lambda\not\in\Lambda(\xb,-\nabla f(\xb);\tilde u)$, we find $\mu\in\Lambda(\xb,-\nabla f(\xb))$ satisfying $\skalp{\mu-\tilde\lambda,\nabla^2g(\xb)(\tilde u,\tilde u)}>0$, and therefore
\begin{equation}\label{EqAux3}
\skalp{\mu-\tilde\lambda,\nabla^2g(\xb)(u_k,u_k)}>0\;\mbox{ for all large }\;k.
\end{equation}
Consider the three possibilities here: $\tilde\lambda=0$, $\tilde\lambda\in\inn\Q^*$, and $\tilde\lambda\in\bd\Q^*\setminus\{0\}$. If $\tilde\lambda=0$, take any sequence $t_k\downarrow 0$ and get by $\tilde u\in\Q$ the equalities
\[
{\rm dist}\big(g(\xb+t_k\tilde u);\Q\big)={\rm dist}\big(g(\xb)+t_k\nabla g(\xb)\tilde u+O(t_k^2);\Q\big)=O(t_k^2)\;\mbox{ as }\;k\to\infty
\]
which allow us, being combined with MSCQ, to find a bounded sequence of $z_k$ satisfying $g(\xb+t_k\tilde u+t_k^2z_k)\in\Q$. Using it together with
$\mu\in\Q^*$ and $\nabla g(\xb)^*\mu=-\nabla f(\xb)=0$ gives us
\begin{equation*}
\begin{array}{ll}
0\ge\skalp{\mu,g(\xb+t_k\tilde u+t_k^2z_k)}&=\disp\Big\la\mu,g(\xb)+\nabla g(\xb)(t_k\tilde u+t_k^2z_k)+\frac{t_k^2}2
\disp\nabla^2g(\xb)(\tilde u,\tilde u)+\oo(t_k^2)\Big\ra\\
&=\disp\Big\la\mu,\frac{t_k^2}2\nabla^2g(\xb)(\tilde u,\tilde u)+\oo(t_k^2)\Big\ra.
\end{array}
\]
Dividing the above inequality by $\frac{t_k^2}2$ and passing to the limit as $k\to\infty$, we conclude that
\[
\skalp{\mu,\nabla^2g(\xb)(\tilde u,\tilde u)}\le 0=\skalp{\tilde\lambda,\nabla^2g(\xb)(\tilde u,\tilde u)}
\]
and thus arrive at a contradiction in the case where $\tilde\lambda=0$.

Assuming now that $\tilde\lambda\in\inn\Q^*$, we get $\lambda^k+\alpha(\mu-\tilde\lambda)\in\inn\Q^*$ for all $\alpha>0$ sufficiently small and all $k\in\N$ sufficiently large. Together with $\nabla g(\xb)^*(\mu-\tilde\lambda)=-\nabla f(\xb)-(-\nabla f(\xb))=0$ it brings us to
$$
\lambda^k+\alpha(\mu-\tilde\lambda)\in\Lambda\big(\xb,x_k^*-\nabla f(\xb)\big)=\Lambda\big(x_k,x_k^*-\nabla f(x_k)\big)\;\mbox{ and}
$$
$$
\skalp{\lambda^k+\alpha(\mu-\tilde\lambda),\nabla^2 g(x_k)(u_k,u_k)}>\skalp{\lambda^k,\nabla^2g(x_k)(u_k,u_k)},
$$
which contradicts the condition $\lambda^k\in\Lambda(x_k,x_k^*-\nabla f(x_k);u_k)$ that follows from the negation of \eqref{EqNeighborTilt}.

Finally, we consider the remaining possibility where $\tilde\lambda\in\bd\Q^*\setminus\{0\}$. Then it follows from the structure of $\Q^*$ that $-\tilde\lambda_0=\norm{\tilde\lambda_r}>0$, $\mu_0+\norm{\mu_r}\le 0$, and
\begin{align*}
\mu_0-\tilde\lambda_0+\disp\Big\la\frac{\tilde\lambda_r}{\norm{\tilde\lambda_r}},\mu_r-\tilde\lambda_r\Big\ra&= \mu_0-\tilde\lambda_0+\Big\la\frac{\tilde\lambda_r}{\norm{\tilde\lambda_r}},\mu_r\Big\ra-\norm{\tilde\lambda_r}\le \mu_0-\tilde\lambda_0+\frac{\norm{\tilde\lambda_r}}{\norm{\tilde\lambda_r}}\norm{\mu_r}-\norm{\tilde\lambda_r}\\
&= \mu_0-\tilde\lambda_0+\norm{\mu_r}-\norm{\tilde\lambda_r}\le 0
\end{align*}
with the usage of the Cauchy-Schwarz inequality. Since the latter holds as equality if and only if the two involved vectors are collinear, we have $\mu_0-\tilde\lambda_0+\skalp{\frac{\tilde\lambda_r}{\norm{\tilde\lambda_r}},\mu_r-\tilde\lambda_r}=0$ if and only if $\mu=\alpha\tilde\lambda$ for some number $\alpha\in\R$. This implies by employing $\nabla g(\xb)^*\mu=\nabla g(\xb)^*\tilde\lambda=-\nabla f(\xb)$ that $(\alpha-1)\nabla g(\xb)^*\tilde\lambda=0$ and consequently that $-\nabla f(\xb)=\nabla g(\xb)^*\tilde\lambda=0$ due to $\mu\not=\tilde\lambda$. The latter contradicts the condition $\tilde\lambda=0$ by $\tilde\lambda\in\Lambda^0(\xb,-\nabla f(\xb))$. Thus we get that
$$
\mu_0-\tilde\lambda_0+\disp\Big\la\frac{\tilde\lambda_r}{\norm{\tilde\lambda_r}},\mu_r-\tilde\lambda_r\Big\ra<0,
$$
which ensures that $\mu_0-\tilde\lambda_0+\skalp{\frac{\lambda_r^k}{\norm{\lambda_r^k}},\mu_r-\tilde\lambda_r}<0$ for all $k$ sufficiently large and hence $\zeta_0^k+\norm{\zeta_r^k}<0$ with $\zeta^k:=\lambda^k+\alpha_k(\mu-\tilde\lambda)$, where $\alpha_k>0$ is chosen to be sufficiently small. Taking into account that $\zeta^k\in\Q^*$ and combining it with $\nabla g(\xb)^*\zeta^k=\nabla g(\xb)^*\lambda^k=x_k^*-\nabla f(\xb)$ tell us that $\zeta^k\in\Lambda(\xb,x_k^*-\nabla f(\xb))$ for all $k\in\N$. It shows together with \eqref{EqAux3} that
\[
\skalp{\zeta^k,\nabla^2g(\xb)(u_k,u_k)}>\skalp{\lambda^k,\nabla^2g(\xb)(u_k,u_k)},\quad k\in\N,
\]
which contradicts the inclusion $\lambda^k\in\Lambda(\xb,x_k^*-\nabla f(\xb);u_k)$ and thus finishes the proof in Case~I.\\[1ex]
{\bf Case~II:} {\em We have  that $x_k\not=\xb$ for all but finitely many $k$}. Suppose without loss of generality that it holds for all $k\in\N$. By passing to a subsequence if necessary, we get that the sequence $(x_k-\xb)/t_k$ converges to some $\bar v\in{\cal S}$, where $t_k:=\norm{x_k-\xb}$. Using
$g(x_k)\in\Q$, $\skalp{\lambda^k,g(x_k)}=0$ for all $k$, and the closedness of the cone $\Q$ leads us to the relationships
\[
\lim_{k\to\infty}\frac{g(x_k)}{t_k}=\lim_{k\to\infty}\frac{g(x_k)-g(\xb)}{t_k}=\nabla g(\xb)\bar v\in\Q\;\mbox{ and }\;0=\lim_{k\to\infty}\frac{\skalp{\lambda^k,g(x_k)-g(\xb)}}{t_k}=\skalp{\tilde\lambda,\nabla g(\xb)\bar v},
\]
which imply that $\bar v\in\K_\Gamma(\xb,-\nabla f(\xb))$ and consequently that $\nabla g(\xb)\bar v=0$ by the assumptions of the theorem.
Denoting $w_k:=t_k^{-2}(x_k-\xb)$ yields
\begin{equation}\label{EqExpansion_g}
g(x_k)=g(\xb)+t_k^2\big(\nabla g(\xb)w_k+\frac 12\nabla^2 g(\xb)(\bar v,\bar v)\big)+\oo(t_k^2)=t_k^2\big(\nabla g(\xb)w_k+\frac 12\nabla^2 g(\xb)(\bar v,\bar v)\big)+\oo(t_k^2).
\end{equation}
To proceed further, we split our analysis into the three subcases:
\begin{equation*}
\limsup_{k\to\infty}\norm{g(x_k)}/t_k^2=\infty,\;\limsup_{k\to\infty}\norm{g(x_k)}/t_k^2=0,\;\mbox{ and }\;0<\limsup_{k\to\infty}\norm{g(x_k)}/t_k^2<\infty,
\end{equation*}
which correspond to the following steps in the proof.\\[1ex]
{\bf Step~1}: {\em Assume that $\limsup_{k\to\infty}\norm{g(x_k)}/t_k^2=\infty$}. Then deduce from \eqref{EqExpansion_g} that $\limsup_{k\to\infty}\norm{\nabla g(\xb)w_k}=\infty$. Passing to a subsequence if necessary, we get that $\lim_{k\to\infty}\norm{\nabla g(\xb)w_k}=\infty$ and that the sequence of $\nabla g(\xb)w_k/\norm{\nabla g(\xb)w_k}$ converges to some element $\nabla g(\xb)w\in{\cal S}$. It follows from \eqref{EqExpansion_g} that
\[
\nabla g(\xb)w=\lim_{k\to\infty}\frac{g(x_k)}{t_k^2\norm{\nabla g(\xb)w_k}}\in\Q\;\mbox{ and }\;\skalp{\tilde\lambda,\nabla g(\xb)w}=\lim_{k\to\infty}\frac{\skalp{\lambda^k,g(x_k)}}{t_k^2\norm{\nabla g(\xb)w_k}}=0,
\]
which yields $w\in\K_\Gamma(\xb,-\nabla f(\xb))$. Together with $\nabla g(\xb)w\not=0$ it contradicts the imposed assumption.\\[1ex]
{\bf Step~2:} {\em Assume that $\limsup_{k\to\infty}\norm{g(x_k)}/t_k^2=0$}. In this setting we deduce from \eqref{EqExpansion_g} that
\[
0=\lim_{k\to\infty}\frac{g(x_k)}{t_k^2}=\lim_{k\to\infty}\nabla g(\xb)w_k+\frac 12\nabla^2g(\xb)(\bar v,\bar v).
\]
Consequently, there is a vector $\tilde w\in\R^n$ such that $\nabla g(\xb)\tilde w+\frac 12\nabla^2g(\xb)(\bar v,\bar v)=0$, and hence the inclusion $(\tilde u,\tilde\lambda,\bar v,\tilde w)\in{\cal Z}$ holds. Next we claim that
\begin{equation}\label{EqAux4}
\lim_{k\to\infty}\lambda_0^k\frac{\norm{\nabla g_r(x_k)u_k}^2-\big(\nabla g_0(x_k)u_k\big)^2}{t_k^2}=0.
\end{equation}
Indeed, assuming $g(x_k)=0$ yields $\nabla g(x_k)u_k\in\Q$, $\lambda^k\in\Q^*$, and $\skalp{\lambda^k,\nabla g(x_k)u_k}=0$ while implying in turn the condition $\lambda_0^k\big(\norm{\nabla g_r(x_k)u_k}^2-(\nabla g_0(x_k)u_k)^2\big)=0$. In the case where $g(x_k)\not=0$ for infinitely many $k$, we deduce from \eqref{EqAux2} that the corresponding sequence of $\skalp{\HH(x_k,\lambda^k)u_k,u_k}$ must be bounded. Then taking into account the definition of $\HH(x_k,\lambda^k)$, Lemma~\ref{LemHH}, and the convergence $g_0(x_k)/t_k^2\to 0$ as $k\to\infty$ verifies the validity of \eqref{EqAux4}.

Our next claim is that condition \eqref{EqAux4} ensures that $\rho(\tilde u,\tilde\lambda,\bar v)=0$. This is certainly true if $\tilde\lambda=0$. In the case where $\tilde\lambda\not=0$ we deduce from $\tilde\lambda_0\le-\norm{\tilde\lambda_r}$ that $\tilde\lambda_0<0$ and therefore
\[
\begin{array}{ll}
0=\disp\lim_{k\to\infty}\frac{\norm{\nabla g_r(x_k)u_k}^2-\big(\nabla g_0(x_k)u_k\big)^2}{t_k^2}=\lim_{k\to\infty}\Big\{&\big\|\nabla
\disp g_r(\xb)\frac{u_k}{t_k}+\nabla^2 g_r(\xb)(\bar v,\tilde u)\big\|^2\\
&-\disp\big(\nabla g_0(\xb)\frac{u_k}{t_k}+\nabla^2 g_0(\xb)(\bar v,\tilde u)\big)^2\Big\}.
\end{array}
\]
Furthermore, we have the equalities
\[
0=\lim_{k\to\infty}\frac{\la\lambda^k,\nabla g(x_k)u_k\ra}{t_k}=\lim_{k\to\infty}\Big\la\tilde\lambda,\nabla g(\xb)\frac{u_k}{t_k}+\nabla^2 g(\xb)(\bar v,\tilde u)\Big\ra.
\]
Supposing that the sequence of $\nabla g(\xb)\frac{u_k}{t_k}$ is unbounded allows us to find a subsequence of $\nabla g(\xb)\frac{u_k}{t_k}/\norm{\nabla g(\xb)\frac{u_k}{t_k}}$, which converges to some $\nabla g(\xb)z\in{\cal S}$ satisfying
$$
\norm{\nabla g_r(\xb)z}^2-\big(\nabla g_0(\xb)z\big)^2=0\;\mbox{ and }\;\skalp{\tilde\lambda,\nabla g(\xb)z}=0.
$$
It follows that, depending on the sign of $\nabla g_0(\xb)z$, either $z\in\K_\Gamma(\xb,-\nabla f(\xb))$ or $-z\in\K_\Gamma(\xb,-\nabla f(\xb))$, which contradicts the assumption above. Thus the sequence of $\nabla g(\xb)\frac{u_k}{t_k}$ is bounded, and so its subsequence converges to some element $\nabla g(\xb)z$ satisfying
$$
\norm{\nabla g_r(\xb)z+\nabla^2 g_r(\xb)(\bar v,\tilde u)}^2-\big(\nabla g_0(\xb)z+\nabla^2 g_0(\xb)(\bar v,\tilde u))^2=0,\;\mbox{ and }\;\skalp{\tilde \lambda,\nabla g(\xb)z+\nabla^2 g(\xb)(\bar v,\tilde u)}=0.
$$
Employing now Lemma~\ref{LemRho} shows that $\rho(\tilde u,\tilde\lambda,\bar v)$ is nonnegative, and therefore we get that $\rho(\tilde u,\tilde\lambda,\bar v)=0$. It follows finally from \eqref{EqAux2} and Lemma~\ref{LemHH} that the inequalities
\[
\bskalp{\big(\nabla^2f(\xb)+\nabla^2\skalp{\tilde\lambda,g}(\xb)\big)\tilde u,\tilde u}\le \limsup_{k\to\infty}\bskalp{\big(\nabla^2f(x_k)+\nabla^2\skalp{\lambda^k,g}(x_k)+\HH(x_k,\lambda^k)\big)u_k,u_k}\le\frac 1\kappa
\]
hold and thus contradict \eqref{EqSecOrder2} in this setting.\\[1ex]
{\bf Step~3:} {\em Assume that $0<\limsup_{k\to\infty}\norm{g(x_k)}/t_k^2<\infty$}. Passing to a subsequence if necessary allows us to use that $g(x_k)/t_k^2$ converges to some vector $q\in\Q\setminus\{0\}$. Then we deduce from \eqref{EqExpansion_g} that $\nabla g(\xb)w_k$ converges to such a vector $\nabla g(\xb)\tilde w$ that $q=\nabla g(\xb)\tilde w+\frac 12\nabla^2g(\xb)(\bar v,\bar v)$. Since $\skalp{\lambda^k,g(x_k)}=0$ for all $k$, it follows that $\skalp{\tilde\lambda,q}=0$ and therefore $(\tilde u,\tilde\lambda,\bar v,\tilde w)\in{\cal Z}$. Then \eqref{EqAux2} together with the nonnegativity of $\skalp{\HH(x_k,\lambda^k)u_k,u_k}$ tells us that the sequence of $\skalp{\HH(x_k,\lambda^k)u_k,u_k}$ is bounded. Let us verify that
\[
\limsup_{k\to\infty}\skalp{\HH(x_k,\lambda^k)u_k,u_k}\ge\frac{\rho(\tilde u,\tilde\lambda,\bar v)}{q_0}.
\]
Indeed, this claim holds trivially if $\tilde\lambda=0$ by the nonnegativity of $\skalp{\HH(x_k,\lambda^k)u_k,u_k}$ and
$\rho(\tilde u,0,\bar v)=0$. Assuming now that $\tilde\lambda\not=0$ yields $\tilde\lambda_0<0$. Taking into account \eqref{EqAux2}, Lemma~\ref{LemHH} as well as the convergence $g_0(x_k)/t_k^2\to q_0>0$ and $\lambda_0^k\to\tilde\lambda_0<0$ tells us that the nonnegative sequence of
\begin{eqnarray*}
\lefteqn{\frac{g_0(x_k)}{-\lambda_0^k t_k^2}\skalp{\HH(x_k,\lambda^k)u_k,u_k}=\frac{\norm{\nabla g_r(x_k)u_k}^2-(\nabla g_0(x_k)u_k)^2}{t_k^2}}\\
&=\disp\Big\|\nabla g_r(\xb)\frac{u_k}{t_k}+\nabla^2 g_r(\xb)(\bar v,\tilde u)+\frac{\oo(t_k)}{t_k}\Big\|^2-\Big(\nabla g_0(\xb)\frac{u_k}{t_k}+\nabla^2 g_0(\xb)(\bar v,\tilde u)+\disp\frac{\oo(t_k)}{t_k}\Big)^2
\end{eqnarray*}
is bounded. The sequence of $\nabla g(x_k)\frac{u_k}{t_k}$ must also be bounded since otherwise a subsequence of
\[
\frac{\nabla g(x_k)u_k/t_k}{\norm{\nabla g(x_k)u_k/t_k}}=\frac{\nabla g(\xb)u_k/t_k+\nabla^2g(\xb)(\bar v,\tilde u)+\oo(t_k)/t_k}{\norm{\nabla g(x_k)u_k/t_k}}
\]
converges to some element $\nabla g(\xb)w\in{\cal S}$ satisfying
\[
\norm{\nabla g_r(\xb)w}^2-\big(\nabla g_0(\xb)w\big)^2=\lim_{k\to\infty}\frac{g_0(x_k)}{-\lambda_0^k t_k^2\norm{\nabla g(x_k)u_k/t_k}^2}\skalp{\HH(x_k,\lambda^k)u_k,u_k}=0\;\mbox{ and}
\]
\[
\skalp{\tilde\lambda,\nabla g(\xb)w}=\lim_{k\to\infty}\frac{\skalp{\lambda^k,\nabla g(x_k)u_k/t_k}}{\norm{\nabla g(x_k)u_k/t_k}}=0.
\]
Thus, depending on the sign of $\nabla g_0(\xb)w$, we have that either $w\in K_\Gamma(\xb,-\nabla f(\xb))$ or $-w\in K_\Gamma(\xb,-\nabla f(\xb))$ with $\nabla g(\xb)w\not=0$, which contradicts the imposed assumption. Hence the sequence of
$$
\nabla g(x_k)\frac{u_k}{t_k}=\nabla g(\xb)\frac{u_k}{t_k}+\nabla^2 g(\xb)(\bar v,\tilde u)+\frac{\oo(t_k)}{t_k}
$$
is bounded, and so is the one of $\nabla g(\xb)u_k/t_k$. Passing to a subsequence if necessary tells us that $\nabla g(\xb)u_k/t_k$ converges to some $\nabla g(\xb)z$. Thus it follows that
\begin{equation*}
\skalp{\tilde\lambda,\nabla g(\xb)z+\nabla^2g(\xb)(\bar v,\tilde u)}=\lim_{k\to\infty}\skalp{\lambda^k,\nabla g(x_k)\frac{u_k}{t_k}}=0.
\end{equation*}
Remembering the definition of $\rho$ in \eqref{EqRho}, we arrive at
\begin{align*}
\limsup_{k\to\infty}\skalp{\HH(x_k,\lambda^k)u_k,u_k}&=\limsup_{k\to\infty}\frac{-\lambda_0^k}{g_0(x_k)/t_k^2}\frac{\norm{\nabla g_r(x_k)u_k}^2-(\nabla g_0(x_k)u_k)^2}{t_k^2}\\
&\ge\frac{-\tilde\lambda_0}{q_0}\Big(\norm{\nabla g_r(\xb)z+\nabla^2g_r(\xb)(\bar v,\tilde u)}^2-\big(\nabla g_0(\xb)z+\nabla^2g_0(\xb)(\bar v,\tilde u)\big)^2\Big)\\
&\ge\frac{\rho(\tilde u,\tilde\lambda,\bar v)}{q_0}.
\end{align*}
Combining it with \eqref{EqAux2} and $q_0=\nabla g_0(\xb)\tilde w+\frac 12\nabla^2g_0(\xb)(\bar v,\bar v)$ yields
\begin{align*}
\lefteqn{\bskalp{\big(\nabla^2f(\xb)+\nabla^2\skalp{\lambda,g}(\xb)\big)\tilde u,\tilde u}+\frac{\rho(\tilde u,\tilde\lambda,\bar v)}{\nabla g_0(\xb)\tilde w+\frac 12\nabla^2g_0(\xb)(\bar v,\bar v)}}\\
&\le\limsup_{k\to\infty}\bskalp{\big(\nabla^2f(x_k)+\nabla^2\skalp{\lambda^k,g}(x_k)+\HH(x_k,\lambda^k)\big)u_k,u_k}\le\frac 1\kappa.
\end{align*}
This contradicts \eqref{EqSecOrder1} and thus completes the proof of the theorem. $\h$\vspace*{0.05in}

The obtained sufficient conditions for tilt stability of local minimizers in SOCPs are rather involved and are hard to verify. Now we derive simplified ones, which imply the sufficient conditions in both Theorems~\ref{ThSuffCondNonZeroCrit} and \ref{ThSuffCondZeroCrit}. The simplified conditions for tilt stability in SOCPs derived in the next theorem are formulated in the same way in the in-kernel and out-of-kernel cases and resemble those established in \cite[Theorem 6.1]{gm15} for NLPs.\vspace*{-0.05in}

\begin{theorem}[\bf simplified sufficient conditions for tilt stability in SOCPs]\label{ThSuffCondSimple} In addition to the assumptions of Theorem~{\rm\ref{ThNeighCharRedSyst}} suppose that, given a number $\kappa>0$, we have the inequality
\[
\bskalp{\big(\nabla^2f(\xb)+\nabla^2\skalp{\lambda,g}(\xb)\big)u,u}>\frac 1\kappa
\]
that is valid for all the multipliers
\[
\lambda\in\tilde\Lambda:=\bigcup_{v\in\K_\Gamma(\xb,-\nabla f(\xb))\cap{\cal S}}\Lambda\big(\xb,-\nabla f(\xb);v\big)
\]
and for all the vectors $u\in{\cal S}$ satisfying
\[
\skalp{\lambda,\nabla g(\xb)u}=0\;\mbox{ and }\;\lambda_0\big(\norm{\nabla g_r(\xb)u}^2-(\nabla g_0(\xb)u)^2\big)=0.
\]
Then $\xb$ is a tilt-stable local minimizer for \eqref{EqProbl} with the prescribed modulus $\kappa$.
\end{theorem}\vspace*{-0.05in}
{\bf Proof.} Let us verify that the conditions imposed in the theorem ensure the fulfillment of the sufficient conditions for tilt stability in both Theorems~\ref{ThSuffCondNonZeroCrit} and \ref{ThSuffCondZeroCrit}. As shown in the proof of Theorem~\ref{ThSuffCondNonZeroCrit}, we have $\Lambda^0(\xb,-\nabla f(\xb))=\{\lb\}$, which implies that $\lb\in\Lambda(\xb,-\nabla f(\xb);\bar u)$ when $\nabla f(\xb)\not=0$. If $\nabla f(\xb)=0$, consider any multiplier $\lambda\in\Lambda(\xb,0)$ and take an arbitrary sequence $t_k\downarrow0$. Since
\[
{\rm dist}\big(g(\xb+t_k\bar u);\Q\big)={\rm dist}\big(t_k\nabla g(\xb)\bar u+\OO(t_k^2);\Q\big)=\OO(t_k^2),
\]
there exists a bounded sequence of $w_k$ satisfying $g(\xb+t_k\bar u+t_k^2w_k)\in\Q$. Taking into account that $\lambda\in\Q^*$ and $\nabla g(\xb)^*\lambda=0$, we get
\begin{equation*}
\begin{array}{ll}
\disp 0\ge\limsup_{k\to\infty}\frac{\skalp{\lambda,g(\xb+t_k\bar u+t_k^2w_k)}}{t_k^2}&=\disp\limsup_{k\to\infty}\frac{\Big\la\lambda,\nabla g(\xb)(t_k\bar u+\disp
t_k^2w_k)+\frac{t_k^2}2\nabla^2 g(\xb)(\bar u,\bar u)\Big\ra}{t_k^2}\\
&=\disp\frac 12\nabla^2\skalp{\lambda,g}(\xb)(\bar u,\bar u).
\end{array}
\end{equation*}
This yields the inclusion $0\in\Lambda(\xb,0;\bar u)$ and thus shows that condition \eqref{EqPointTilt1} of Theorem~\ref{ThSuffCondNonZeroCrit} holds.

Next we verify that both conditions (a) and (b) of Theorem~\ref{ThSuffCondZeroCrit} are fulfilled. The validity of (a) follows from the fact that
for every $u\in\K_\Gamma(\xb,-\nabla f(\xb))$ and every $\lambda\in\Lambda(\xb,-\nabla f(\xb);u)$ we have $\skalp{\lambda,\nabla g(\xb)u}=0$ and $\lambda_0\big(\norm{\nabla g_r(\xb)u}^2-(\nabla g_0(\xb)u)^2\big)=0$. To verify (b), consider an arbitrary quadruple $(u,\lambda,v,w)\in{\cal Z}$ and observe that for any $\mu\in\Lambda(\xb,-\nabla f(\xb))$ we get
\[
\Big\la\mu-\lambda,\frac 12\nabla^2 g(\xb)(v,v)\Big\ra=\Big\la\mu-\lambda,\nabla g(\xb)w+\frac 12\nabla^2 g(\xb)(v,v)\Big\ra=\Big\la\mu,\nabla g(\xb)w+\frac 12\nabla^2 g(\xb)(v,v)\Big\ra\le 0,
\]
which yields $\lambda\in\Lambda(\xb,-\nabla f(\xb);v)$. Then (b) follows by taking into account that $\rho(u,\lambda,v)\ge 0$ by Lemma~\ref{LemRho} and $\nabla g_0(\xb)w+\frac 12\nabla^2g_0(\xb)(v,v)>0$ due to $0\not=\nabla g(\xb)w+\frac 12\nabla^2g(\xb)(v,v)\in\Q$. $\h$\vspace*{0.07in}

It is important to observe that the {\em nondegeneracy} condition, which means that $\nabla g(\xb)$ has full rank, implies that the multiplier set $\Lambda(\xb,-\nabla f(\xb))$ is a {\em singleton} and that the term $\rho(u,\lambda,v)$ {\em vanishes}. Furthermore, in this case the sufficient conditions obtained in Theorems~\ref{ThSuffCondNonZeroCrit}, \ref{ThSuffCondZeroCrit}, and \ref{ThSuffCondSimple} are {\em equivalent}. The next example demonstrates that all these phenomena {\em fail} for programs that exhibit {\em degeneracy}.\vspace*{-0.05in}

\begin{example}[\bf tilt stability under degeneracy]\label{ex1} {\rm Consider the following program of type \eqref{EqProbl}:
\begin{eqnarray*}
\begin{array}{ll}
\mbox{minimize }\;f(x):=\disp\frac{1}{4}\left(3 x_1^2+\frac{7-\sqrt{15}}{\sqrt{15}}\big(x_2^2+x_1x_2\big)\right)-x_3\;\mbox{ with }\;x=(x_1,x_2,x_3)\in\R^3\\
\mbox{subject to}\;g(x)=\big(g_1(x),g_2(x),g_3(x)\big)\in\Q,
\end{array}
\end{eqnarray*}
where the functions $g_i(x)$, $i=1,2,3$, are defined by
\[
g_1(x):=\frac{1}{2}(x_1^2+x_2^2+x_1x_2),\quad g_2(x):=\frac{1}{4}\left(\frac{1}{2}x_1^2+x_2^2+x_1x_2\right),\quad g_3(x):=\frac{1}{4}\big(x_1^2+x_2^2+x_1x_2\big)+x_3.
\]
It is not hard to check by direct calculations that $\ox:=(0,0,0)$ is a stationary point of the SOCP under consideration. We intend to show that Theorem~\ref{ThSuffCondSimple} cannot detect that $\ox$ is a tilt-stable local minimizer of this problem while the more complicated Theorem~\ref{ThSuffCondZeroCrit} can.

Let us first check that MSCQ is fulfilled for the constraint system $g(x)\in\Q$ at $\xb$ represented by
$$
\Q=\big\{q=(q_0,q_r)\in\R\times\R^2\;\big|\;h(q):=\norm{q_r}-q_0\le 0\big\}.
$$
For any $x_1,x_2$ denote $\alpha:=g_1(x_1,x_2,x_3)$ and $\beta:=g_2(x_1,x_2,x_3)$ and observe that $h(g(x_1,x_2,\hat x_3))=0$ for $\hat x_3:=-\alpha/2\pm \sqrt{\alpha^2-\beta^2}$. We have furthermore that $x=(x_1,x_2,x_3)\in g^{-1}(\Q)$ if and only if
\[
-\frac{\alpha}{2}-\sqrt{\alpha^2-\beta^2}\le x_3\le-\frac{\alpha}{2}+\sqrt{\alpha^2-\beta^2}> 0.
\]
Consider now $x_3:=-\alpha/2+\sqrt{\alpha^2-\beta^2}+\gamma$ with some $\gamma>0$ and note that $g_3(x_1,x_2,x_3)=\sqrt{\alpha^2-\beta^2}+\gamma$. Taking into account that $\sqrt{3}\alpha\le 2\sqrt{\alpha^2-\beta^2}$ and $\beta\le\alpha/2$ and denoting by $L>0$ a Lipschitz constant of $h$, we get the distance estimate
\begin{eqnarray*}
\dist{\left(x;g^{-1}(\Q)\right)}&\le&\gamma=\frac{2}{\sqrt{3}}\left(\sqrt{\left(\alpha+\frac{\sqrt{3}}{2}\gamma\right)^2}-\alpha\right)
\le\frac{2}{\sqrt{3}}\left(\sqrt{\alpha^2+2\sqrt{\alpha^2-\beta^2}\gamma+\gamma^2}-\alpha\right)\\
&=&\frac{2}{\sqrt{3}}\left(\sqrt{\beta^2+\left(\sqrt{\alpha^2-\beta^2}+\gamma\right)^2}-\alpha\right)
=\frac{2}{\sqrt{3}}\big|h\big(g(x)\big)\big|\le\frac{2}{\sqrt{3}}L\,\dist{\left(g(x);\Q\right)}.
\end{eqnarray*}
The case where $x_3:=-\alpha/2-\sqrt{\alpha^2-\beta^2}+\gamma$ with some $\gamma<0$ can be treated similarly, and hence the claimed metric subregularity is verified.

To proceed further, calculate the needed values $\nabla f(\xb)=(0,0,-1)$,
\[
\K_\Gamma\big(\xb,-\nabla f(\xb)\big)=\ker\nabla g(\xb)=\R^2\times\{0\},\;\mbox{ and }\;\Lambda\big(\xb,-\nabla f(\xb)\big)=\big\{(a,b,1)\;\big|\;a\le-\sqrt{b^2+1}\big\}.
\]
Thus, given a direction $(t,s,0)\in\K_\Gamma(\xb,-\nabla f(\xb))$, we get that $(a,b,1)\in\Lambda(\xb,-\nabla f(\xb);(t,s,0))$ if and only if the pair
$(a,b)$ is a solution to the optimization problem:
\[
\mbox{maximize }\;a(t^2+s^2+ts)+\frac{b}{2}\left(\frac{t^2}{2}+s^2+ts\right)+\frac{1}{2}\big(t^2+s^2+ts\big)\;\mbox{ subject to }\;a\le-\sqrt{b^2+1}.
\]
Taking into account that $t^2+s^2+ts>0$, we can reduce the latter problem to the unconstrained maximization with respect to $b$ only:
\[
\mbox{maximize }\;-\sqrt{b^2+1}\big(t^2+s^2+ts\big)+\frac{b}{2}\left(\frac{t^2}{2}+s^2+ts\right)+\frac{1}{2}\big(t^2+s^2+ts\big),
\]
which can be solved explicitly. Denoting $A:=t^2+s^2+ts$ and $B:=t^2/2+s^2+ts$, the optimal value of $b$ and consequently of $a$ are calculated by
\begin{equation}\label{eq : ab_def}
b=\frac{B}{\sqrt{4A^2-B^2}}\;\mbox{ and }\;a=-\frac{2A}{\sqrt{4A^2-B^2}}.
\end{equation}
Note that the numbers $A$, $B$, and $A-B$ are nonnegative, and hence it follows that $b\ge 0$. Furthermore, for any multiplier $\lambda=(a,b,1)$ and any direction $u=(p,q,0)\in\ker\nabla g(\xb)$ we have
\begin{equation}\label{eq : SecOrd}
\bskalp{\big(\nabla^2f(\xb)+\nabla^2\skalp{\lambda,g}(\xb)\big)u,u}=\left(a+\frac{b}{4}+2\right)p^2+\left(a+\frac{b}{2}+\frac{7}{2\sqrt{15}}\right)\big(q^2+pq\big).
\end{equation}
Let us check that the number in \eqref{eq : SecOrd} is nonnegative for all $\lambda\in\Lambda(\xb,-\nabla f(\xb);v)$ with some $v=(t,s,0)\in\K_\Gamma(\xb,-\nabla f(\xb))\cap{\cal S}$. Indeed, note that \eqref{eq : ab_def} yields $(a+b/4+2)>0$ and $(a+b/4+2)>(a+b/2+7/2\sqrt{15})$ by $B^2\le A^2$. Since $t^2+s^2=1$, we get that $A=\mu B$ for some $\mu\in[1,2]$, and so
\[
\sqrt{4A^2-B^2}\left(a+\frac{b}{2}+\frac{7}{2\sqrt{15}}\right)=\left(\frac{7}{2\sqrt{15}}\sqrt{4\mu^2-1}+\frac{1}{2}-2\mu\right)B\ge 0
\]
for all $\mu\in[8/11,2]\supset[1,2]$. The equality is attained therein if and only if $\mu=2$, which gives us the multiplier $\tilde\lambda:=(-4/\sqrt{15},1/\sqrt{15},1)$ and the corresponding directions equal either $v=(\pm 1,0,0)$ or $v=(\pm\sqrt{2}/2,\mp\sqrt{2}/2,0)$. This ensures therefore that for every multiplier $\lambda\ne\tilde\lambda=(-4/\sqrt{15},1/\sqrt{15},1)$ the number in \eqref{eq : SecOrd} is strictly positive whenever direction $u=(p,q,0)$ with $p^2+q^2=1$ is taken. Indeed, it follows from the estimate
\[
\begin{array}{ll}
\disp\left(a+\frac{b}{4}+2\right)p^2+\left(a+\frac{b}{2}+\frac{7}{2\sqrt{15}}\right)\disp\big(q^2+pq\big)&\ge\disp\left(a+\frac{b}{2}+\frac{7}{2\sqrt{15}}\right)
\disp\big(p^2+q^2+pq\big)\\
&\ge\disp\frac{1}{2}\left(a+\frac{b}{2}+\frac{7}{2\sqrt{15}}\right).
\end{array}
\]
Considering now the noted multiplier $\tilde\lambda=(a,b,1)=(-4/\sqrt{15},1/\sqrt{15},1)$, we have that $(a+b/2+7/2\sqrt{15})=0$, and so the number in
\eqref{eq : SecOrd} equals $(a+b/4+2)p^2\ge 0$. Choosing the directions $u=(0,\pm 1,0)$ gives us the equality
\[
\bskalp{\big(\nabla^2f(\xb)+\nabla^2\skalp{\tilde\lambda,g}(\xb)\big)u,u}=0,
\]
which shows that the tilt stability of $\xb$ cannot be detected by Theorem~\ref{ThSuffCondSimple}.

On the other hand, with the same $\tilde\lm$ and $u=(0,\pm 1,0)$ we get $\rho(u,\tilde\lambda,v)>0$ for all the directions
$v=(k,l,0)$ such that $k+2l\ne 0$ and $k^2+l^2=1$. Indeed, using the values $\skalp{\nabla^2g(\xb)u,v}=\pm(k/2+l)(1,1/2,1/2)$
and $\nabla g(\xb)z=(0,0,z_3)$ for any $z=(z_1,z_2,z_3)$ and taking into account that
$\skalp{\tilde\lambda,\nabla g(\xb)z+\skalp{\nabla^2g(\xb)u,v}}=0$ tell us that $z_3=\pm(7-\sqrt{15})/2\sqrt{15}(k/2+l)$. Consequently, it yields
\[
\nabla g(\xb)z+\skalp{\nabla^2g(\xb)u,v}=\pm\left(\frac{k}{2}+l\right)\left(1,\frac{1}{2},\frac{7}{2\sqrt{15}}\right)
\]
and $49/60(k/2+l)^2+(k/2+l)^2/4-(k/2+l)^2>0$ if $k/2+l\ne 0$. This verifies that $\rho(u,\tilde\lambda,v)>0$.

It follows from the proof of Theorem~\ref{ThSuffCondSimple} that $\lambda\in\Lambda(\xb,-\nabla f(\xb);v)$ for any quadruple $(u,\lambda,v,w)\in{\cal Z}$.
Since it is the case for the chosen vectors $\tilde\lambda$, $u=(0,\pm 1,0)$, and $v$ equal to either $(\pm 1,0,0)$ or $v=(\pm \sqrt{2}/2,\mp\sqrt{2}/2,0)$, with some $w\in\R^n$, we get $\rho(u,\tilde\lambda,v)>0$ and thus deduce from Theorem~\ref{ThSuffCondZeroCrit} in the in-kernel case that $\ox$ is a tilt-stable minimizer for the SOCP under consideration.}
\end{example}\vspace*{-0.25in}

\section{Pointbased Necessary Conditions and Criteria for Tilt Stability}\vspace*{-0.05in}

This section is mainly devoted to deriving pointbased necessary conditions for tilt-stable minimizers of SOCPs. The obtained results are complementary to the sufficient conditions for such minimizers given in Section~4 and, being unified with the latter, allow us to establish {\em complete pointbased characterizations} of tilt stability in second-order cone programming.

As in Section~4, it makes sense to consider separately the in-kernel and out-of-kernel cases. We start with the easier out-of-kernel case in which the following {\em no-gap} necessary condition holds while being different from the sufficient one in Theorem~\ref{ThSuffCondNonZeroCrit} only by the nonstrict inequality sign.\vspace*{-0.05in}

\begin{theorem}[\bf necessary condition for tilt-stable minimizers of SOCPs in the out-of-kernel case]\label{ThNecCondNonZeroCrit} Let
$\xb$ be a tilt-stable local minimizer of \eqref{EqProbl} with modulus $\kappa$ under the standing assumptions of Theorem~{\rm\ref{ThNeighCharRedSyst}}. Suppose in addition that there is a vector $\ub\in\K_\Gamma(\xb,-\nabla f(\xb))$ satisfying $\nabla g(\xb)\ub\not=0$ and take the multiplier $\bar\lambda$ defined in \eqref{barl}. Then for every direction $u\in{\cal S}$ satisfying $\skalp{\lb,\nabla g(\xb)u}=0$ and $\lb_0(\norm{\nabla g_r(\xb)u}^2-(\nabla g_0(\xb)u)^2)=0$ we have the condition
\begin{equation}\label{EqNecPointTilt1}
\bskalp{\big(\nabla^2 f(\xb)+\nabla^2\skalp{\lb,g}(\xb)\big)u,u}\ge\frac 1\kappa.
\end{equation}
Furthermore, the exact bound of tilt stability of \eqref{EqProbl} at $\ox$ is lower estimated by
\begin{align*}
{\rm tilt}(f+\delta_\Gamma,\xb)\ge\sup\Bigl\{1\big/\bskalp{\big(\nabla^2f(\xb)+\nabla^2\skalp{\lb,g}(\xb)\big)u,u}\;\Big|\;&u\in{\cal S},\;\skalp{\lb,\nabla g(\xb)u}=0,\\
&\lb_0\big(\norm{\nabla g_r(\xb)u}^2-(\nabla g_0(\xb)u)^2\big)=0\Bigr\}.
\end{align*}
\end{theorem}\vspace*{-0.05in}
{\bf Proof.} We proceed with the verification of the necessary condition \eqref{EqNecPointTilt1} while observing that the exact bound estimate follows directly from the proof below. Assuming on the contrary that condition \eqref{EqNecPointTilt1} fails for the tilt-stable minimizer $\ox$, find $\tilde u\in{\cal S}$ satisfying
$$
\skalp{\lb,\nabla g(\xb)\tilde u}=0,\;\lb_0\big(\norm{\nabla g_r(\xb)u}^2-(\nabla g_0(\xb)\tilde u)^2\big)=0,\;\mbox{ and }\;\bskalp{\big(\nabla^2f(\xb)+\nabla^2\skalp{\lb,g}(\xb)\big)\tilde u,\tilde u}<\frac 1\kappa
$$
and then show that $\xb$ is not a tilt-stable local minimizer of \eqref{ThNeighCharRedSyst} with modulus $\kappa$ by using the neighborhood characterization of this property taken from Theorem~\ref{ThNeighCharRedSyst}. We proceed with considering the following two possible settings.

Suppose first that $\nabla f(\xb)\not=0$, and so $\nabla g(\xb)^*\lb=-\nabla f(\xb)\not=0$ with $\lb\not=0$. This yields
$$
\nabla g(\xb)\ub\in\bd\Q\setminus\{0\},\;\norm{\nabla g_r(\xb)\ub}=\nabla g_0(\xb)\ub>0,\;\mbox{ and }\;-\lb_0=\norm{\lb_r}>0.
$$
Therefore $\vert\nabla g_0(\xb)\tilde u\vert=\norm{\nabla g_r(\xb)\tilde u}$ and, depending on the sign of $\nabla g_0(\xb)\tilde u$, we have that either $\tilde u\in\K_\Gamma(x,-\nabla f(\xb))$ or $-\tilde u\in\K_\Gamma(x,-\nabla f(\xb))$. Since the left-hand side of \eqref{EqNecPointTilt1} is quadratic in $u$, assume without loss of generality that $\tilde u\in\K_\Gamma(x,-\nabla f(\xb))$. Using the same arguments as in the proof of Theorem~\ref{ThSuffCondNonZeroCrit} leads us to the equalities
$$
\Lambda^0\big(\xb,-\nabla f(\xb)\big)=\Lambda\big(\xb,-\nabla f(\xb)\big)=\big\{\lb\big\},
$$
which imply that $\lb\in\Lambda(\xb,-\nabla f(\xb);\tilde u)$ and show that \eqref{EqNeighborTilt} is violated with $u=\tilde u$, $x=\xb$, $x^*=0$, and $\lambda=\lb$ for every $\eta>0$. This verifies that $\xb$ is not a tilt-stable local minimizer of \eqref{EqProbl} with modulus $\kappa$.

Now suppose that $\nabla f(\xb)=0$, and hence $\lb=0$ by \eqref{barl}. Choosing an arbitrary sequence $t_k\downarrow 0$ and using MSCQ give us
$x_k\in\Gamma$ satisfying the estimate
\[
\norm{x_k-(\xb+t_k\bar u)}\le\sigma{\rm dist}\big(g(\xb+t_k\ub);\Q\big)=\oo(t_k)\;\mbox{ for all large }\;k\in\N.
\]
If $g(x_k)\in\inn\Q$ holds for infinitely many $k$, then for these $k$ we have $\tilde u\in\K_\Gamma(x_k,0)=\R^n$ and $\bskalp{\big(\nabla^2 f(x_k)+\nabla^2\skalp{0,g}(x_k)\big)\tilde u,\tilde u}<1/\kappa$ when $k$ is sufficiently large. Thus the neighborhood characterization \eqref{EqNeighborTilt} for the tilt stability of $\xb$ with modulus $\kappa$ fails whenever $\eta>0$ with $u=\tilde u$, $x=x_k$, $x^*=\nabla f(x_k)$, $\lambda=0$, and large $k$. Suppose further that $g(x_k)\in\bd\Q$ for all but finitely many $k$. Since
$$
g(x_k)=g(\xb)+\nabla g(\xb)(x_k-\xb)+\oo(\norm{x_k-\xb})=t_k\nabla g(\xb)\ub+\oo(t_k),
$$
it follows that $g_0(x_k)=\norm{g_r(x_k)}>0$ for all large $k$. Next we consider the real-valued function $h(x):=\norm{g_r(x)}-g_0(x)$ and observe that the condition
$$
\nabla h(x_k)=\frac{g_r(x_k)^*}{\norm{g_r(x_k)}}\nabla g_r(x_k)-\nabla g_0(x_k)\not=0
$$
yields the existence of $x_k'$ which is arbitrary close to $x_k$ and such that $h(x_k')<0$, i.e., $g(x_k')\in\inn\Q$. Thus we are in the same position as just before that gives us a contradiction with \eqref{EqNeighborTilt}.

The remaining case is where $\nabla h(x_k)=0$. By the assumed MSCQ we get from \eqref{EqTanGamma} that $T_\Gamma(x_k)=\{u\in\R^n\mid\nabla g(x_k)u\in T_\Q(g(x_k))\}$, which together with
$$
T_\Q(g(x_k))=\Big\{v\in\R^n\;\Big|\;\Big\la\frac{g_r(x_k)}{\norm{g_r(x_k)}},v_r\Big\ra-v_0\le 0\Big\}\supset\Big(-1,\frac{g_r(x_k)}{\norm{g_r(x_k)}}\Big)^\perp
$$
ensures the following relationships:
\[
T_\Gamma(x_k)\supset\Big\{u\in\R^n\;\Big|\;\Big\la\frac{g_r(x_k)}{\norm{g_r(x_k)}},\nabla g_r(x_k)u\Big\ra-\nabla g_0(x_k)u=0\Big\}=\big\{u\in\R^n\;\big|\;\nabla h(x_k)u=0\big\}=\R^n.
\]
Thus $\tilde u\in\K_\Gamma(x_k,0)$, and for every $\eta>0$ condition \eqref{EqNeighborTilt} is violated again with $u=\tilde u$, $x=x_k$, $x^*=\nabla f(x_k)$, and $\lambda=0$ when $k$ is chosen sufficiently large. This completes the proof of the theorem. $\h$\vspace*{0.05in}

Combining the results of Theorems~\ref{ThSuffCondNonZeroCrit} and \ref{ThNecCondNonZeroCrit}, we arrive at the following effective pointbased characterization of tilt stability for SOCPs in the out-of-kernel case.\vspace*{-0.05in}

\begin{theorem}[\bf pointbased characterization of tilt-stable minimizers for SOCPs in the out-of-{\nobreak kernel} case]\label{char-out} In the setting of Theorem~{\rm\ref{ThNeighCharRedSyst}}, assume that there exists $\ub\in\K_\Gamma(\xb,-\nabla f(\xb))$ such that $\nabla g(\xb)\ub\not=0$. Then $\xb$ is a tilt-stable local minimizer for \eqref{EqProbl} with some modulus $\kappa>0$ if and only if for every direction $u\in{\cal S}$ satisfying $\skalp{\lb,\nabla g(\xb)u}=0$ and $\lb_0(\norm{\nabla g_r(\xb)u}^2-(\nabla g_0(\xb)u)^2)=0$ we have
\begin{equation*}
\bskalp{\big(\nabla^2 f(\xb)+\nabla^2\skalp{\lb,g}(\xb)\big)u,u}>0,
\end{equation*}
where the multiplier $\lb$ is defined in \eqref{barl}. Furthermore, the exact bound for tilt stability of \eqref{EqProbl} at $\ox$ is precisely calculated by the formula
\begin{align*}
{\rm tilt\,}(f+\delta_\Gamma,\xb)=\sup\Bigl\{1\big/\bskalp{\big(\nabla^2 f(\xb)+\nabla^2\skalp{\lb,g}(\xb)\big)u,u}\;\Big|\;&u\in{\cal S},\;\skalp{\lb,\nabla g(\xb)u}=0,\\ &\lb_0\big(\norm{\nabla g_r(\xb)u}^2-(\nabla g_0(\xb)u)^2\big)=0\Bigr\}.
\end{align*}
\end{theorem}\vspace*{-0.05in}
{\bf Proof.} The major observation here is that the modulus of tilt stability is not specified in the formulation of the theorem, which concerns therefore tilt-stable local minimizers of  \eqref{EqProbl} with {\em some} modulus $\kk>0$. Thus the claimed characterization and exact bound formula are derived from combining the corresponding results of Theorems~\ref{ThSuffCondNonZeroCrit} and \ref{ThNecCondNonZeroCrit}. $\h$\vspace*{0.05in}

Prior to examining the in-kernel case in what follows, we present the following technical lemma that concerns a certain error bound for the second-order cone constraint in \eqref{EqProbl}.\vspace*{-0.05in}

\begin{lemma}[\bf error bound for second-order cone constraint systems]\label{Lem_h} Define $h\colon\R^n\to\R$ by $h(x):=\norm{g_r(x)}-g_0(x)$. Given any $\delta>0$, there is a number $\epsilon>0$ and a neighborhood $U$ of $\xb$ such that for every $x_0\in U$ with $|h(x_0)|<\epsilon\norm{g_r(x_0)}$ and $\norm{\nabla h(x_0)}\ge\delta$ we can find $\tilde x$ satisfying
$$
h(\tilde x)=0\;\mbox{ and }\norm{\tilde x-x_0}\le 2\frac{\vert h(x_0)\vert}\delta.
$$
\end{lemma}\vspace*{-0.1in}
{\bf Proof.} Let $L>0$ be a common Lipschitz constant for the mappings $g(x)$ and $h(x)$ on $\xb+\B$. We can see that $h$ is $C^2$-smooth around every $x$ with $g_r(x)\not=0$ and its derivatives are calculated by
\begin{align*}
\nabla h(x)&=\frac{g_r(x)^*}{\norm{g_r(x)}}\nabla g_r(x)-\nabla g_0(x),\\
\nabla^2 h(x)&=\frac{g_r(x)^*}{\norm{g_r(x)}}\nabla^2g_r(x)-\nabla^2 g_0(x)+\frac 1{\norm{g_r(x)}}\nabla g_r(x)^*\big(I-\frac{g_r(x)g_r(x)^*}{\norm{g_r(x)}^2}\big)\nabla g_r(x),
\end{align*}
where $I$ is the identity matrix. Thus we can find constants $A,B>0$ such that
\[
\norm{\nabla^2 h(x)}\le A+\frac{B}{\norm{g_r(x)}}\;\mbox{ for all }\;x\in\ox+\B\;\mbox{ with }\;g_r(x)\not=0.
\]
Take an arbitrary number $\delta>0$ and then choose $\gg>0$ and $0<\epsilon\le\delta/4L$ so small that
$$
\frac{A\vert h(x)\vert+2B\epsilon}{\delta^2}\le\frac 12\;\mbox{ and }\;\gg+2\frac{\vert h(x)\vert}\delta\le 1\;\mbox{ for all }\;x\in\B_\gg(\ox).
$$
Pick $x_0\in\B_\gg(\ox)$ such that $\vert h(x_0)\vert<\epsilon\norm{g_r(x_0)}$ and $\norm{\nabla h(x_0)}\ge\delta$ and denote $\rho:=2\vert h(x_0)\vert/\delta$. Then there exists $x_1$ satisfying $h(x_0)+\nabla h(x_0)(x_1-x_0)=0$ and $\norm{x_1-x_0}\le\vert h(x_0)\vert/\delta$. Furthermore, for all $x\in\B_\rho(x_0)\subset\ox+\B$ we get the inequalities
\[
\begin{array}{ll}
\disp\norm{g_r(x)}\ge\norm{g_r(x_0)}-L\norm{x-x_0}&\ge\norm{g_r(x_0)}-L\rho=\disp\norm{g_r(x_0)}-2L\frac{\vert h(x_0)\vert}\delta\\
&\ge\disp\norm{g_r(x_0)}\Big(1-2L\frac\epsilon\delta\Big)\ge\frac{\norm{g_r(x_0)}}2,
\end{array}
\]
which bring us to the estimate
\[
\norm{\nabla^2h(x)}\le A+\disp\frac{2B}{\norm{g_r(x_0)}}.
\]
It implies that $\nabla h(\cdot)$ is Lipschitz continuous on $\B_\rho(x_0)$ with constant $A+2B/\norm{g_r(x_0)}$ and that
\[
\tau:=\frac 1{\delta}\Big(A+\frac{2B}{\norm{g_r(x_0)}}\Big)\norm{x_1-x_0}\le\Big(A+\frac{2B}{\norm{g_r(x_0)}}\Big)\frac{\vert h(x_0)\vert}{\delta^2}\le\frac{A\vert h(x_0)\vert+2B\epsilon}{\delta^2}\le\frac 12.
\]
Denoting $\hat\tau:=\tau^{-1}(1-\sqrt{1-2\tau})\le 2$, we apply Robinson's stability result from \cite[Theorem~2]{Rob72} to get the existence of $\tilde x\in x_0 +\hat\tau\vert h(x_0)\vert/\delta\B\subset\B_\rho(x_0)$ with $h(\tilde x)=0$. This completes the proof. $\h$\vspace*{0.05in}

To proceed further with deriving necessary conditions for and then characterizations of tilt-stable minimizers in the in-kernel case of SOCPs, we invoke the notion of {\em 2-regularity}, which was initiated by Tret'yakov \cite{Tr84} in the case of zero Jacobian and then was strongly developed by Avakov \cite{Av85} whose definition is mainly reflected in what follows.\vspace*{-0.05in}

\begin{definition}[\bf 2-regularity of mappings]\label{Def2Regularity}Let $g:\R^n\to\R^s$ be twice Fr\'echet differentiable at $\xb\in\R^n$. We say that
$g$ is {\sc 2-regular} at the point $\xb$ in the {\sc direction} $v\in\R^n$ if for any $p\in\R^s$ the system
\[
\nabla g(\xb)z+\nabla^2g(\xb)(v,w)=p,\;\nabla g(\xb)w=0
\]
admits a solution $(z,w)\in\R^n\times\R^n$.
\end{definition}
\vspace*{-0.05in}

The next theorem provides two necessary conditions for tilt-stable minimizers in SOCPs with a prescribed modulus in the in-kernel setting. Both of them are parallel to the sufficient conditions of Theorem~\ref{ThSuffCondZeroCrit} with the replacement of the strict inequality therein by the nonstrict one. Note that only the second condition below requires the assumption on 2-regularity in addition to the standing MSCQ.\vspace*{-0.05in}

\begin{theorem}[\bf necessary condition for tilt-stable minimizers of SOCPs in the in-kernel case]\label{ThNecCondZeroCrit} Let $\xb$ be a tilt-stable local minimizer for \eqref{EqProbl} with modulus $\kappa$ under the standing assumptions of Theorem~{\rm\ref{ThNeighCharRedSyst}}, and assume that for every critical direction $u\in\K_\Gamma(\xb,-\nabla f(\xb))$ we have $\nabla g(\xb)u=0$. Suppose further that for every $(u,\lambda,v,w)\in{\cal Z}$ fulfilling $\nabla g(\xb)w+\frac 12\nabla^2 g(\xb)(v,v)=0$ the mapping $g(\cdot)$ is 2-regular in the direction $v$. Then the following conditions are satisfied:

{\bf (a)} For every $u\in K_\Gamma(\xb,-\nabla f(\xb))\cap{\cal S}$ and every $\lambda\in\Lambda(\xb,-\nabla f(\xb);u)$ we have
\[
\bskalp{\big(\nabla^2f(\xb)+\nabla^2\skalp{\lambda,g}(\xb)\big)u,u}\ge\frac 1\kappa.
\]

{\bf(b)} For every $(u,\lambda,v,w)\in{\cal Z}$ 
 we have
\begin{equation}\label{EqNecSecOrder1}
\begin{array}{ll}
\disp\bskalp{\big(\nabla^2f(\xb)+\nabla^2\skalp{\lambda,g}(\xb)\big)u,u}+\frac{\rho(u,\lambda,v)}{\nabla g_0(\xb)w+\frac 12\nabla^2 g_0(\xb)(v,v)}\ge\frac 1\kappa\\
\disp\mbox{ whenever }\;\nabla g(\xb)w+\frac 12\nabla^2 g(\xb)(v,v)\not=0;\quad\mbox{ and}
\end{array}
\end{equation}\vspace*{-0.15in}
\begin{equation}\label{EqNecSecOrder2}
\disp\bskalp{\big(\nabla^2f(\xb)+\nabla^2\skalp{\lambda,g}(\xb)\big)u,u}\ge\frac 1\kappa\;\mbox{ whenever }\;\rho(u,\lambda,v)=0.
\end{equation}
Furthermore, the lower estimate of the exact bound of tilt stability is given by
$$
{\rm tilt}(f+\delta_\Gamma,\xb)\ge\frac 1{\min\{\chi_1,\chi_2,\chi_3\}},
$$
where the numbers $\chi_i$, $i=1,2,3$, are taken from Theorem~{\rm\ref{ThSuffCondZeroCrit}}.
\end{theorem}\vspace*{-0.05in}
{\bf Proof.} Note first that the claimed lower estimate of the tilt stability bound follows from the proofs of the necessity of conditions (a) and (b). To verify the necessity of (a) for tilt stability, we just pass to the limit as $\eta\dn 0$ in the neighborhood characterization (ii) of Theorem~\ref{ThNeighCharRedSyst}. To prove the necessity of (b), pick a quadruple $(u,\lambda,v,w)\in{\cal Z}$ and get for any $t>0$ that
\[
{\rm dist}\big(g(\xb+tv+t^2w);\Q\big)={\rm dist}\Big(g(\xb)+t\nabla g(\xb)v+t^2\big(\nabla g(\xb)w+\frac 12\nabla^2g(\xb)(v,v)\big)+\oo(t^2);\Q\Big)=\oo(t^2).
\]
We split the subsequent proof into the following two cases.\\[1ex]
{\bf Case~I:} $\lambda=0$. In this case $\nabla f(\xb)=0$, and by the assumed MSCQ for every $t>0$ we find some $x_t\in\Gamma$ satisfying $\norm{x_t-(\xb+tv+t^2w)}=\oo(t^2)$. If there is a sequence $t_k\downarrow 0$ such that $g(x_{t_k})\not=0$, then it follows from the proof of Theorem~\ref{ThNecCondNonZeroCrit} that $\skalp{\nabla^2 f(\xb)u,u}\ge 1/\kappa$. This verifies \eqref{EqNecSecOrder2} since $\rho(u,0,v)=0$. Suppose now that
$g(x_t)=0$ for all small $t>0$. By
$$
g(x_t)=g\big(\xb+tv+t^2w+\oo(t^2)\big)=t^2\Big(\nabla g(\xb)w+\disp\frac 12\nabla^2g(\xb)(v,v)\Big)+\oo(t^2)
$$
we have $\nabla g(\xb)w+\frac 12\nabla^2 g(\xb)(v,v)=0$, and thus $g(\cdot)$ is 2-regular at $\ox$ in the direction $v$. Fixing any $q\in\inn\Q$ and applying \cite[Proposition~2(c)]{GfrOut16} ensure the existence of $\beta>0$ such that
\[
{\rm dist}\big(x_t;g^{-1}(t^3q)\big)\le\frac\beta{\norm{x_t-\xb}}\norm{t^3q-g(x_t)}=\OO(t^2)\;\mbox{ for small }\;t>0.
\]
Thus we find points $x_t'$ for which $\norm{x_t'-x_t}=\OO(t^2)$ and $g(x_t')=t^3q\in\inn\Q$. This tells us that $\K_\Gamma(x_t',0)=\R^n$. Using further
Theorem~\ref{ThNeighCharRedSyst} with $x=x_t'$, $x^*=\nabla f(x_t')$, and $\lambda=0$ yields $\skalp{\nabla^2 f(x_t')u,u}\ge\frac 1\kappa$ for all small $t>0$.  To get finally \eqref{EqNecSecOrder2} in this case, we pass to the limit as $t\dn 0$ with taking into account the equality $\rho(u,0,v)=0$.\\[1ex]
{\bf Case~II:} $\lambda\not=0$. Since $\lambda\in\Q^*$, it follows that $-\lambda_0\ge\norm{\lambda_r}$, and consequently $-\lambda_0>0$. We have furthermore that $\vert\nabla g_0(\xb)u\vert=\norm{\nabla g_r(\xb)u}$ and $\skalp{\lambda,\nabla g(\xb)u}=0$, which implies that either $u\in\K_\Gamma(\xb,-\nabla f(\xb))$ or $-u\in\K_\Gamma(\xb,-\nabla f(\xb))$. This gives us $\nabla g(\xb)u=0$ by the assumption of the theorem. We split the subsequent analysis in this case into the following two steps.\\[1ex]
{\bf Step~1:} $\nabla g(\xb)w+\frac 12\nabla^2g(\xb)(v,v)=0$. Then $g(\xb+tv+t^2w)=\oo(t^2)$, and its follows from \cite[Proposition~2(c)]{GfrOut16} due to the imposed 2-regularity of $g(\cdot)$ at $\ox$ in the direction $v$ that there is $\beta>0$ for which
\[
{\rm dist}\big(\xb+tv+t^2w;g^{-1}(0)\big)\le\frac{\beta}{\norm{tv+t^2w}}\norm{g(\xb+tv+t^2w)}=\oo(t)\;\mbox{ for small }\;t>0.
\]
This allows us to find for such $t$ some vectors $x_t\in\R^n$ for which $g(x_t)=0$ and $\norm{x_t-(\xb+tv+t^2w)}=\oo(t)$. If $\rho(u,\lambda,v)\not=0$, then there is nothing to prove, and hence we suppose that $\rho(u,\lambda,v)=0$. Employing Lemma~\ref{LemRho} in this case gives us a vector $z\in\R^n$ satisfying
\[
-\lambda_0\Big(\norm{\nabla g_r(\xb)z+\nabla^2g_r(\xb)(v,u)}^2-\big(\nabla g_0(\xb)z+\nabla^2g_0(\xb)(v,u)\big)^2\Big)=\skalp{\lambda,\nabla g(\xb)z+\nabla^2g(\xb)(v,u)}=0,
\]
which implies in turn the equalities
\[
\nabla g(x_t)(u+tz)=\big(\nabla g(\xb)+t\nabla ^2g(\xb)v+\oo(t)\big)(u+tz)=\nabla g(\xb)u+t\big(\nabla g(\xb)z+\nabla^2g(\xb)(v,u)\big)+\oo(t)=\oo(t).
\]
It follows from \cite[Proposition~2(d)]{GfrOut16} that there is a number $\beta'>0$ such that for all $t>0$ sufficiently small we get vectors $u_t\in\R^n$ satisfying
\[
\nabla g(x_t)u_t=0\;\mbox{ and }\;\norm{u_t-(u+tz)}\le\beta'\frac{\norm{\nabla g(x_t)(u+tz)}}{\norm{x_t-\xb}}=\frac{\oo(t)}t,
\]
which tells us that $u_t\in\K_\Gamma(x_t,\nabla g(x_t)^*\lambda)$. Applying \cite[Proposition~2(d)]{GfrOut16} again ensures that the Jacobian matrix $\nabla g(x_t)$ is of full rank, and therefore
$$
\Lambda\big(x_t,\nabla g(x_t)^*\lambda\big)=\Lambda\big(x_t,\nabla g(x_t)^*\lambda;u_t\big)=\big\{\lambda\big\}.
$$
Employing the neighborhood condition \eqref{EqNeighborTilt1} from Theorem~\ref{ThNeighCharRedSyst} with $x=x_t$, $x^*=\nabla f(x_t)+\nabla g(x_t)^*\lambda$, and $u=u_t/\norm{u_t}$ and then taking into account that $\HH(x_t,\lambda)=0$ due to $g(x_t)=0$, we arrive at
$$
\Big\la\big(\nabla^2f(x_t)+\nabla^2\skalp{\lambda,g}(x_t)\big)\frac{u_t}{\norm{u_t}},\frac{u_t}{\norm{u_t}}\Big\ra\ge\frac 1\kappa\;\mbox{ for all small}\;t>0.
$$
Passing there to the limit as $t\dn0$ verifies \eqref{EqNecSecOrder2} in this setting.\\[1ex]
{\bf Step~2:} $\nabla g(\xb)w+\frac 12\nabla^2g(\xb)(v,v)\not=0$. Remembering that $\lm\ne 0$ and $\lambda\in N_\Q(\nabla g(\xb)w+\frac 12\nabla^2 g(\xb)(v,v))$ in this case gives us the conditions
$$
\lambda\in\bd\Q^*,\;\nabla g(\xb)w+\frac 12\nabla^2g(\xb)(v,v)\in\bd\Q,\;\mbox{ and }\;\hat\lambda=\alpha\big(\nabla g(\xb)w+\frac 12\nabla^2 g(\xb)(v,v)\big)
$$
for some $\alpha>0$. It implies, in particular, that $-\lambda_0=\norm{\lambda_r}>0$ and $\nabla g_0(\xb)w+\frac 12\nabla^2 g_0(\xb)(v,v)=\norm{\nabla g_r(\xb)w+\frac 12\nabla^2 g_r(\xb)(v,v)}>0$. Note also that $\nabla g(\xb)^*\lambda\not=0$, since otherwise $\nabla f(\xb)=0$ and consequently $\lambda=0$. For the function $h(x)=\norm{g_r(x)}-g_0(x)$ we have the representations
\begin{align*}
h(\xb+tv+t^2w)&=\norm{g_r(\xb+tv+t^2w)}- g_0(\xb+tv+t^2w)\\
&=\norm{g_r(\xb)+t\nabla g_r(\xb)v+t^2(\nabla g_r(\xb)w+\frac 12\nabla^2 g_r(\xb)(v,v))}\\
&\quad-\Big(g_0(\xb)+t\nabla g_0(\xb)v+t^2\big(\nabla g_0(\xb)w+\frac 12\nabla^2 g_0(\xb)(v,v)\big)\Big)+\oo(t^2)=\oo(t^2),
\end{align*}\vspace*{-0.15in}
\begin{align*}
\nabla h(\xb+tv+t^2w)&=\frac{g_r(\xb+tv+t^2w)^*}{\norm{g_r(\xb+tv+t^2w)}}\nabla g_r(\xb+tv+t^2w)-\nabla g_0(\xb+tv+t^2w)\\
&=\frac{\big(t^2(\nabla g_r(\xb)w+\frac 12\nabla g_r(\xb)(v,v))+\oo(t^2)\big)^*}{\norm{t^2\big(\nabla g_r(\xb)w+\frac 12\nabla g_r(\xb)(v,v)\big)+\oo(t^2)}}\big(\nabla g_r(\xb)+\OO(t)\big)-\big(\nabla g_0(\xb)+\OO(t)\big)\\
&=\disp\frac{\big(\lambda_r+\disp\oo(t^2)/t^2\big)^*}{\norm{\lambda_r}+\disp\oo(t^2)/t^2}\big(\nabla g_r(\xb)+\OO(t)\big)-\big(\nabla g_0(\xb)+\OO(t)\big)=\frac{\lambda^*}{\norm{\lambda_r}}\nabla g(\xb)+\frac{\oo(t^2)}{t^2},
\end{align*}
which yield, in particular, the relationships
$$
\lim_{t\downarrow 0}\norm{\nabla h(\xb+tv+t^2w)}=\frac{\norm{\lambda^*\nabla g(\xb)}}{\norm{\lambda_r}}=\frac{\norm{\nabla f(\xb)}}{\norm{\lambda_r}}>0\;\mbox{ and }\;\lim_{t\downarrow 0}\norm{\nabla h(x_t)}=\frac{\norm{\lambda^*\nabla g(\xb)}}{\norm{\lambda_r}}>0.
$$
Using now Lemma~\ref{Lem_h} gives us vectors $x_t$ satisfying $h(x_t)=0$ and $\norm{x_t-(\xb+tv+t^2w)}=\oo(t^2)$ for small $t>0$. Defining for such $t$ the multipliers
$\lambda^t:=\alpha\hat g(x_t)/t^2\in N_\Q(g(x_t))$ and employing
\[
g(x_t)=g(\xb+tv+t^2w)+\oo(t^2)=t^2\Big(\nabla g(\xb)w+\frac 12\nabla^2 g(\xb)(v,v)\Big)+\oo(t^2)
\]
ensure that $\lim_{t\downarrow 0}\lambda^t=\lambda$. Since for $\rho(u,\lambda,v)=\infty$ condition \eqref{EqNecSecOrder1} certainly holds, we consider the case where $\rho(u,\lambda,v)$ is finite. Then Lemma~\ref{LemRho} gives us $z\in\R^n$ satisfying the equalities
$$
\skalp{\lambda,\nabla g(\xb)z+\nabla^2g(\xb)(v,u)}=0,\;\rho(u,\lambda,v)=-\lambda_0\Big(\norm{\nabla g_r(\xb)z+\nabla^2g_r(\xb)(v,u)}^2-\big(\nabla g_0(\xb)z+\nabla^2g_0(\xb)(v,u)\big)^2\Big).
$$
Remembering that $g_0(x_t)=\norm{g_r(x_t)}$ leads us to the expressions
\begin{align*}
\nabla h(x_t)(u+tz)&=\frac{g_r(x_t)^*}{\norm{g_r(x_t)}}\nabla g_r(x_t)(u+tz)-\nabla g_0(x_t)(u+tz)=\Big\la\frac{\hat g(x_t)}{\norm{g_r(x_t)}},\nabla g(x_t)(u+tz)\Big\ra\\
&=\Big\la\frac{\hat g(x_t)}{\norm{g_r(x_t)}},\nabla g(\xb)u+t\big(\nabla g(\xb)z+\nabla^2g(\xb)(v,u)\big)+\oo(t)\Big\ra\\
&=t\Big\la\frac{\lambda^t}{\norm{\lambda^t_r}},\nabla g(\xb)z+\nabla^2g(\xb)(v,u)\Big\ra+\oo(t)\\
&=t\Big\la\frac{\lambda^t}{\norm{\lambda^t_r}}-\frac{\lambda}{\norm{\lambda_r}},\nabla g(\xb)z+\nabla^2g(\xb)(v,u)\Big\ra+\oo(t)=\oo(t).
\end{align*}
In this way we get, whenever $t>0$ is sufficiently small, that
$$
\nabla h(x_t)u_t=\frac 1{\norm{\lambda^t_r}}\big\la\lambda^t,\nabla g(x_t)u_t\big\ra=0\quad\mbox{for}\quad u_t:=u+tz-\frac{\nabla h(x_t)(u+tz)}{\norm{\nabla h(x_t)}^2}\nabla h(x_t)^*=u+tz+\oo(t).
$$
It follows that $u_t\in\K_\Gamma(x_t,x_t^\ast-\nabla f(x_t))$, where $x_t^\ast:=\nabla f(x_t)+\nabla g(x_t)^*\lambda^t\to 0$ as $t\to 0$. Since $N_\Q(g(x_t))=\{\alpha\hat g(x_t)\mid\alpha\ge 0\}$ and $x_t^\ast-\nabla f(x_t)\not=0$, we easily see that $\Lambda(x_t,x_t^\ast-\nabla f(x_t))=\{\lambda^t\}$. Therefore, the neighborhood condition \eqref{EqNeighborTilt1} from Theorem~\ref{ThNeighCharRedSyst} implies that
\begin{equation}\label{EqAux5}
\Big\la\Big(\nabla^2f(x_t)+\nabla^2\skalp{\lambda^t,g}(x_t)+\HH(x_t,\lambda^t)\Big)\frac{u_t}{\norm{u_t}},\frac{u_t}{\norm{u_t}}\Big\ra\ge\frac 1\kappa\;\mbox{ for small }\;t>0.
\end{equation}
Using the construction of the curvature function $\HH$ in \eqref{H} together with the limiting relations $\lambda^t_0\to\lambda_0$ and $g_0(x_t)/t^2\to\nabla g_0(\xb)w+\frac 12\nabla^2g_0(\xb)(v,v)$ as $t\downarrow 0$, we get
\begin{align*}
\lim_{t\downarrow 0}&\bskalp{\HH(x_t,\lambda^t)u_t,u_t}=\lim_{t\downarrow 0}\frac{-\lambda^t_0}{g_0(x_t)}\Big(\norm{\nabla g_r(x_t)u_t}^2-(\nabla g_0(x_t)u_t)^2\Big)\\
&=\lim_{t\downarrow 0}\frac{-\lambda^t_0}{g_0(x_t)}\Big(\norm{\nabla g_r(\xb)u+t\big(\nabla g_r(\xb)z+\nabla^2g_r(\xb)(v,u)\big)+\oo(t)}^2\\
&\qquad\qquad-\big(\nabla g_0(\xb)u +t\big(\nabla g_0(\xb)z+\nabla^2g_0(\xb)(v,u)\big)+\oo(t)\big)^2\Big)\\
&=\lim_{t\downarrow 0}\frac{-\lambda^t_0}{g_0(x_t)/t^2}\Big(\norm{\nabla g_r(\xb)z+\nabla^2g_r(\xb)(v,u)+\frac{\oo(t)}t}^2-\big(\nabla g_0(\xb)z+\nabla^2g_0(\xb)(v,u)+\frac{\oo(t)}t\big)^2\Big)\\
&=\frac{\rho(u,\lambda,v)}{\nabla g_0(\xb)w+\frac 12\nabla^2 g_0(\xb)(v,v)}.
\end{align*}
Finally, we arrive at \eqref{EqNecSecOrder1} by passing to the limit in \eqref{EqAux5} as $t\dn 0$. $\h$\vspace*{0.05in}

We conclude this section by deriving the following verifiable characterization of tilt-stable minimizers for SOCPs in the in-kernel case.\vspace*{-0.05in}

\begin{theorem}[\bf pointbased characterization of tilt-stable minimizers for SOCPs in the in-kernel case]\label{char-in} In addition to the assumptions of Theorem~{\rm\ref{ThNeighCharRedSyst}} suppose that $\nabla g(\xb)u=0$ for every critical direction $u\in\K_\Gamma(\xb,-\nabla f(\xb))$ and that for every quadruple $(u,\lambda,v,w)\in{\cal Z}$ with $\nabla g(\xb)w+\frac 12\nabla^2 g(\xb)(v,v)=0$ the mapping $g(\cdot)$ is 2-regular in the direction $v$. Then $\xb$ is a tilt-stable local minimizer for problem \eqref{EqProbl} with some modulus $\kk>0$ if and only if the following conditions hold simultaneously:

{\bf (a)} For every $u\in K_\Gamma\big(\xb,-\nabla f(\xb)\big)\cap{\cal S}$ and every $\lambda\in\Lambda(\xb,-\nabla f(\xb);u)$ we get
\[
\bskalp{\big(\nabla^2f(\xb)+\nabla^2\skalp{\lambda,g}(\xb)\big)u,u}>0.
\]

{\bf(b)} For every quadruple $(u,\lambda,v,w)\in{\cal Z}$ we get
\begin{equation*}
\begin{array}{ll}
\disp\bskalp{\big(\nabla^2f(\xb)+\nabla^2\skalp{\lambda,g}(\xb)\big)u,u}+\frac{\rho(u,\lambda,v)}{\nabla g_0(\xb)w+\frac 12\nabla^2 g_0(\xb)(v,v)}>0\\
\disp\mbox{ whenever }\;\nabla g(\xb)w+\frac 12\nabla^2 g(\xb)(v,v)\not=0;\quad\mbox{ and}
\end{array}
\end{equation*}\vspace*{-0.15in}
\begin{equation*}
\disp\bskalp{\big(\nabla^2f(\xb)+\nabla^2\skalp{\lambda,g}(\xb)\big)u,u}>0\;\mbox{ whenever }\;\rho(u,\lambda,v)=0.
\end{equation*}
Furthermore, the exact bound of tilt-stability of $\ox$ in \eqref{EqProbl} is calculated by the formula
\[
{\rm tilt}(f+\delta_\Gamma,\xb)=\frac 1{\min\{\chi_1,\chi_2,\chi_3\}},
\]
where the numbers $\chi_i$, $i=1,2,3$, are taken from Theorem~{\rm\ref{ThSuffCondZeroCrit}}.
\end{theorem}\vspace*{-0.05in}
{\bf Proof.} We derive the claimed statements from the results of Theorems~\ref{ThSuffCondZeroCrit} and \ref{ThNecCondZeroCrit} by taking into account that modulus of tilt stability is not specified in this theorem. $\h$\vspace*{-0.15in}

\section{Concluding Remarks}\vspace*{-0.05in}

The results of this paper and the presented proofs show that the situation with deriving pointwise conditions (sufficient, necessary, and characterizations) for tilt-stable minimizers in nonpolyhedral cone programming without nondegeneracy (and thus with nonunique Lagrange multiplies) is a highly challenging and dramatically more involved task in comparison with polyhedral and/or nondegenerate settings. Our attention to this issue is motivated not only by the theoretical interest but also by the increasing importance of tilt stability for the design and justification of primal-dual numerical algorithms of optimization and obtaining their convergence rates.

Based on advanced recent developments in second-order variational analysis, we are able to establish here comprehensive qualitative and quantitative results on tilt stability in SOCPs generated by a single second-order cone under merely metric subregularity constraint qualification that is far removed not only from nondegeneracy but also from the conventional Robinson constraint qualification in conic programming. However, it remains to study SOCPs described by products of second-order cones and to investigate other remarkable classes of conic programs. Note that the recent results of \cite{GfrMo17} allows us to calculate the main generalized differential construction of second-order variational analysis implemented here, the so-called subgradient graphical derivative, for the normal cone mappings generated by constraint sets in $C^2$-reducible conic programs. This makes it possible to proceed with deriving neighborhood characterizations of tilt-stable minimizers in such programs based on the abstract criterion in \cite{ChHiNg17} and the technique developed above. On the other hand, establishing pointbased conditions for tilt stability required by applications is a big issue for our future research.\vspace*{0.05in}

{\bf Acknowledgements.} The research of the first and second authors was partially supported by the Austrian Science Fund (FWF) under grant P29190-N32. The research of the third author was partially supported by the USA National Science Foundation under grants DMS-1512846 and DMS-1808978, and by the USA Air Force Office of Scientific Research under grant No.\,15RT0462.

\end{document}